\documentclass[12pt]{article}
\usepackage{amssymb}
\begin{document}
\title{\LARGE\bf 
A triality group 
of nonassociative algebras \\
with involution \\
 \ \\
\large{In memorry of Professor Susumu Okubo (1930-2015)} \\}
\date{}
\author{\large \bf Noriaki Kamiya$^1$ and Susumu Okubo$^2$\\ \\
$^1$Department of Mathematics, 
University of Aizu, \\
Aizuwakamatsu, Japan\\
$^2$Department of Physics and Astronomy, University of Rochester,\\
Rochester, N.Y,  U.S.A}
\maketitle
\thispagestyle{empty}
\vskip 2mm
{\bf Abstract}
\par
In this paper, we give concept of triality groups and study
a characterization of 
the group for
symmetric composition a;gebras.
\vskip 3mm
AMS classification (2010), 17A30, 17B40, 20F29.
\par
Keywords, Symmetric composition algebras, triality groups, local and global triality relation.
\par 
\vskip 5mm
{\bf Introduction}
\par
The main purpose of this study is to exhibit a generalization of the automorphism group with respect to nonassociative algebras
with involution.
Moreover frankly describing,
we will provide concept of new idea or frontier  of a triality relation.
for a group called a global triality relation.
\par
First
let $A$ be a nonassociative algebra over a field $F$.
If a triple $\sigma=(\sigma_{1},\sigma_{2},\sigma_{3})
\in (Epi (A))^{3}$, where $Epi (A)$ is the set of onto and endomorphisms of $A$,
satisfies
$\sigma_{j}(xy)=(\sigma_{j+1}x)(\sigma_{j+2}y)$
for
$j=1,2,3$
with
$\sigma_{j\pm 3}=\sigma_{j}$
and for any $x,y\in A,$
then set of all such
triples forms
a group called the triality group of $A$ in componentwise multiplication
and is denoted by $Trig(A)$.
In particular, we consider a construction of the group for
 symmetric composition algebras $A$
and try  to find  sub-groups of $Trig (A)$.
As a bi-product of left and right operators of $A$,
we have found a formula
for invariant sub-groups of  automorphism
groups of symmetric composition algebras ( to see, Theorem 4.1 and Corollary 4.2). 
\par
Next summarizing this note, we will discuss several sections as in
the following.
\par
1. Triality groups (definitions and preliminary)
\par
2. Symmetric composition algebras
\par
3. Local triality relations
\par
4. Automorphisms of symmetric composition algebras
\par
5, Examples of triality groups for some nonassociative algebras 
with involution
\par
Appendix. simple examples
\par
This article may be viewed as a continuation
of the previous paper one ([K-O.15]).
However the content is independent of previous work,
that is, we will describe in order to readers as self-contained as possible.  
\par
Finally, we would like to note that 
this work is one in the end of  his life (S.Okubo, 1930-2015).
\par
\vskip 3mm
{\bf 1\ Triality Group (Definitions and Preliminary)}
\par
\vskip 3mm
Let $A$ be an algebra over a field $F$ of charachteristic
not $2$
with a bi-linear product denoted by
juxtaposition
$xy$ for
$x,y\in A $.
Suppose that a triple $g=(g_{1},g_{2},g_{3})\in (Epi A)^{3}$,
where $Epi  A$ denote the set of epimorphisms of $A$,
satisfies a global triality relation
$$
g_{j}(xy)=
(g_{j+1}x)(g_{j+2}y)\eqno(1.1)$$
for any
$x,y\in A$ and for any
$j=1,2,3$,
Here the index $j$
is defined by modulo $3$
so that
$$
g_{j\pm3}=g_{j}.\eqno(1.2)
$$
For the second triple
$g^{'}=(g_{1}^{'},g_{2}^{'},g_{3}^{'})\in (Ep{i} A)^{3}$
satisfying the same triality relation,
we introduce their product componetwise by
$$
gg^{'}=
(g_{1}g^{'}_{1},g_{2}g^{'}_{2},g_{3}g^{'}_{3}).
\eqno(1.3)$$
Then,
a set consisting of all such 
triples forms a group which we call
the triality group of $A$,
and write
$$
{ Trig}(A)=
\{
g=(g_{1},g_{2},g_{3})\in
(Ep{i}A)^{3}|
g_{j}(xy)=
(g_{j+1}x)
(g_{j+2}y),
\forall x,y\in A,\ \forall \ j=1,2,3\}.
\eqno(1.4)$$
Here we emphasize that this group
 is clearly a generalization of the 
automorphism group defined by
$$
{Auto}(A)=
\{g\in Ep{i}A|
g(xy)=
(gx)(gy),
\forall x,y\in A\}.\eqno(1.5)$$
We then note that ${ Trig}(A)$
is invariant under actions of an alternative
group $A_{4}$
(or equivalently the tetrahedral group )
as follows.
First,
let 
$\phi\in End ({\rm Trig}(A))$
by
$$
\phi:g_{1}\rightarrow g_{2}\rightarrow g_{3}\rightarrow g_{1}
\eqno(1.6)
$$
which satisfies
$\phi^{3}=id$
and leaves Eq.(1.1)
invariant.
Thus,
${ Trig}(A)$
is invariant under actions of the cyclic group $Z_{3}$
generated by $\phi$.
We next introduce
$\tau_{\mu}\in End({ Trig}(A))$
for $\mu=1,2,3$
by
$$
\tau_{1}:g_{1}\rightarrow g_{1},\ g_{2}\rightarrow -g_{2},\ 
g_{3}\rightarrow -g_{3}$$
$$
\tau_{2}:g_{1}\rightarrow -g_{1},\ g_{2}\rightarrow g_{2},\ 
g_{3}\rightarrow -g_{3}
\eqno(1.7)$$
$$
\tau_{3}:g_{1}\rightarrow -g_{1},\ g_{2}\rightarrow -g_{2},\ 
g_{3}\rightarrow g_{3},
$$
which leave
Eq.(1.1)
invariant again.
Moreover,
they satisfy
$$
\tau_{\mu}\tau_{\nu}=\tau_{\nu}\tau_{\mu},\ 
\tau^{2}_{\mu}=id,\ 
\tau_{1}\tau_{2}\tau_{3}=id\eqno(1.8)
$$
for 
$\mu,\nu=1,2,3,$
so that the group 
 generated by $<id,\tau_{1},\tau_{2},\tau_{3}>$
is isomorphic to the Klein's
$4$-group $K_{4}.$
\par
Further we note
$$
\phi\tau_{\mu}\phi^{-1}=
\tau_{\mu+1},\ 
({\rm with}\ \tau_{4}=\tau_{1}).\eqno(1.9)$$
Since
$Z_{3}$
and
$K_{4}$
generate the alternative group
$A_{4}$,
this shows that
$A_{4}$
is a invariant group of
$Trig (A)$, for example, 
$(\tau_{1}g_{1},\tau_{1}g_{2},\tau_{1}g_{3})=(g_{1},-g_{2},-g_{3})
\in Trig (A).$
\par
If $A$ is involutive with the involution map
$x\rightarrow {\bar x}$
satisfying
$$
\bar{\bar x}=x,\ \overline{ xy}=
{\bar y}{\bar x},\eqno(1.10)$$
we define
${\bar Q}\in End(A)$
for any
$Q\in End(A)$
by
$$
\overline{Qx}=
{\bar Q}{\bar x}.\eqno(1.11)
$$
then,
taking the involution of Eq.(1.10)
and letting
$x\leftrightarrow {\bar y}$,
we find
$$
{\bar g}_{j}(xy)=
({\bar g}_{j+2}x)
({\bar g}_{j+1}y),\eqno(1.12)
$$
so that
$\theta\in End({ Trig}(A))$
defined by
$$
\theta:
g_{1}\rightarrow{\bar g}_{2},\ g_{2}
\rightarrow {\bar g}_{1},\ 
g_{3}\rightarrow{\bar g}_{3}\eqno(1.13)
$$
yields also a invariant operation of 
$Trig (A)$.
Moreover,
we obtain
$$
\phi\theta\phi=\theta,\ 
\theta^{2}=id,\ 
\theta\tau_{1}\theta^{-1}=\tau_{2},\ 
\theta\tau_{2}\theta^{-1}=\tau_{1},\ 
\theta\tau_{3}\theta^{-1}=\tau_{3}.\eqno(1.14)
$$
Then,
$A_{4}$
togehter with
$\theta$
give the
$S_{4}$-symmetry for
$Trig(A)$
with identification  of
$$
\phi=(1,2,3),\ 
\tau_{1}=(2,3)(1,4)\ 
\tau_{2}=(3.1)(2.4),
\tau_{3}=(1,2)(3,4),
\theta=(1,2)\eqno(1.15)
$$
in terms of the transpositions of the
$S_{4}$-group.
This also shows that
$S_{4}$ is a invaruant group of $Trig(A)$, if $A$ is involutive. 
\par
Before going into further discussion,
we note that
$Trig (A)$ for any
$A$ always contains a Klein's 4-group $K_{4}$.
In fact, 
let $Id\in Epi(A)$
be defined by
$(Id)x=x,$
for $x\in A,$
and set
$$
{\tilde \tau}_{0}=
(Id,Id,Id),\ 
{\tilde \tau}_{1}=(Id.-Id,-Id),\ 
{\tilde \tau}_{2}=
(-Id,Id,-Id),\ 
{\tilde\tau}_{3}=
(-Id,-Id,Id).\eqno(1.16)$$
Then,
these $4$ element of $(Epi(A)^{3})$
satisfy Eq.(1.8)
in view of 
Eq.(1.3),
so that they give another Klein's $4$-group.
Moreover,
they satisfy Eq.(1.1).
Similarly,
we have ${ Trig}(F)=K_{4},$
for the simplest
case of $A=F.$
\par
In contrast to the global triality relation Eq.(1.1),
we may also consider the local triality relation
$$
t_{j}(xy)=(t_{j+1}x)y+x(t_{j+2}y)\eqno(1.17)
$$
for $t_{j}\in End(A)$
with
$t_{j\pm3}=t_{j}.$
Analogously to
Eq.(1.4),
we introduce
$$
s\circ Lrt(A)=
$$
$$
\{
t=(t_{1},t_{2},t_{3})\in
(End\ A)^{3}|
t_{j}(xy)=
(t_{j+1}x)y+
x(t_{j+2}y),
\forall x,y\in A,
\forall i=1,2,3\}.\eqno(1.18)$$
Then,
it defines a Lie algebra now with
component-wise commutation relation
(]K-O.15]).
\par
Here,
$s\circ\ Lrt(A)$
stands for symmetric Lie-related triple,
which has been reformed to as ${\rm stri}(A)$
in
[O.05]
instead.
To see this symmetric Lie- related triple, if we set
$$
t_{j}^{'}=
\sum^{3}_{k=1}\alpha_{j-k}t_{k},\ (j=1,2,3)\eqno(1.19)$$
for any
$\alpha_{j}\in F$ satisfying
$\alpha_{j\pm3}=\alpha_{j},$
then it is easy to verify that we have also
$$
t^{'}=(t_{1}^{'},t_{2}^{'},t_{3}^{'})\in s\circ Lrt(A).\eqno(1.20)
$$
If the exponential map $t_{j} \rightarrow \xi_{j}$
is given by
$$
\xi_{j}={\rm exp}(t_{j})=
\sum^{\infty}_{n=0}
{1\over n!}
(t_{j})^{n}
,
\eqno(1.21)
$$
is well defined,
then we have shown in
[K-O.15]
that the validity of
$$
\xi_{j}(xy)=
(\xi_{j+1}x)(\xi_{j+2}y),\eqno(1.22)$$
provided that
$t=(t_{1},t_{2},t_{3})\in s\circ Lrt(A)$
and vice-versa.
\par
Note that the existence ot the exponantial map
requires the underlying field
$F$ to be at least of zero charachteristic.
\par
We first introduce multiplication 
operators
$L(x),R(x)\in End\ A$
by
$$
L(x)y=xy,\ R(x)y=yx\eqno(1.23)$$
as usual.
Then,
Eq.(1.17) yields
$$
t_{j}L(x)=L(x)
t_{j+2}+L(t_{j+1}x),\eqno(1.24a)$$
$$
t_{j}R(y)=R(y)
t_{j+1}+R(t_{j+2}y)\eqno(1.24b)$$
while Eq.(1.1)
gives
$$
g_{j}L(x)=L(g_{j+1}x)g_{j+2},\eqno(1.25a)$$
$$
g_{j}R(y)=R(g_{j+2}y)g_{j+1}\eqno(1.25b)$$
any $g=(g_{1},g_{2},g_{3})\in { Trig}(A).$
\par
From Eqs.(1.25),
we find
$$
g_{j}L(x)R(y)g_{j}^{-1}=
L(g_{j+1}x)R(g_{j+1}y)\eqno(1.26a)$$
$$
g_{j}R(y)L(x)g_{j}^{-1}=R(g_{j+2}y)L(g_{j+2}x).\eqno(1.26b)$$
We next introduce the notion of a regular triality algebra.
\par
\vskip 3mm
{\bf Def.1.1}
\par
\vskip 3mm
Let $d_{j}(x,y)\in End(A)$
for
$x, y\in A$ and for
$j=1,2,3$ be
to satisfy
\par
(i)
$$
d_{1}(x,y)=R(y)L(x)-R(x)L(y),\eqno(1.27a)$$
$$
d_{2}(x,y)=
L(y)R(x)-L(x)R(y).\eqno(1.27b)$$
\par
(ii)\quad
The explicit form for
$d_{3}(x,y)$
is unspecified except for
$$
d_{3}(y,x)=-d_{3}(x,y)\eqno(1.27c)$$
\par
(iii)\quad
$(d_[{}(x,y),d_{2}(x,y),d_{3}(x,y))\in s\circ Lrt(A),$
i.e.,
they satisfy
$$
d_{j}(x,y)(uv)=(d_{j+1}(x,y)u)v+
u(d_{j+2}(x,y)v)\eqno(1.28)
$$
for any
$x,y,u,v\in A$
and for any
$j=1,2,3.$
here the index
over $j$
is defined modulo $3$ as before.
\par
We call the algebra $A$ 
satisfying these conditions to be a regular
triality algebra.
([K-O.15])
\par
For the case of Lie algebra $A$  equipped with the product $[x,y]$ as usual, if we set
$$
d_{1}(x,y)=R(y)L(x)-R(x)L(y),$$
$$
d_{2}(x,y)=L(y)R(x)-L(x)R(y),$$
$$
d_{3}(x,y)=d_{0}(x.y)=L([x,y]),$$
where $L(x)y=[x,y]$ and $R(x)y=[y,x]$,
then we get
$$d_{j}=L([x,y]) {\rm \  for\  any }\  j=0,1,2.
$$
Thus this implies that $d_{j}\   ^{'}s$ is a inner derivation of $A$. 
( i.e., a special example of the Def.1.1). 
\par
\vskip 3mm
{\bf Def.1.2}
\par
\vskip 3mm
{\bf Condition (B)}:
We have
$AA=A.$
\par
\vskip 3mm
{\bf Condition (C)}:
If some $b\in A$
satisfies either
$L(b)=0,$
or
$R(b)=0,$
then
$b=0.$
\par
We can now prove.
\par
\vskip 3mm
{\bf Proposition 1.3}
\par
\vskip 3mm
\it
Let $A$ be a regular triality
algebra
satisfying either the
condition (B)
or (C).
We then obtain the followings:
\par
(i)\quad
 For any $t=(t_{1},t_{2},t_{3})\in s\circ Lrt(A)$,
we have
$$
[t_{j},d_{k}(x,y)]=
d_{k}(t_{j-k}x,y)+
d_{k}(x,t_{j-k}y).\eqno(1.29a)$$
Especially,
if we choose
$t_{j}=d_{j}(u,v),$
it yields also
$$
[d_{j}(u,v),d_{k}(x,y)]=
d_{k}(d_{j-k}(u,v)x,y)+
d_{k}(x,d_{j-k}(u,v)y).\eqno(1.29b)$$
\par
(ii)\quad
For any $g=(g_{1},g_{2},g_{3})\in { Trig}(A),$
we have
$$
g_{j}d_{k}(x,y)g_{j}^{-1}=
d_{k}(g_{j-k}x, g_{j-k}y)\eqno(1.30)
$$
for any
$j,k=1,2,3$
and for any
$u,v,x,y\in A,$
\par
\vskip 3mm
\rm
{\bf Proof}
\par
\vskip 3mm
Since Eqs.(1.29) have been already proved in
[K-O,15],
we will give only a proof of Eq.(1.30)
below.
In view of Eq.(1.26)
and
(1.27),
we see that Eq.(1.30)
holds valid for any
$j=1,2,3$
and for
$k=1,2.$
Therefore,
it suffices
to prove of the case of $k=3.$
To this end,
we set
$$
D_{j,k}:\equiv
g_{j}d_{k}(x,y)g_{j}^{-1}-d_{k}(g_{j-k}x,g_{j-k}y)\eqno(1.31)$$
for a fixed
$x,y\in A.$
Then,
as we have noted,
we have
$D_{j,1}=D_{j,2}=0$
identically.
Moreover,
we will show that it satisfies also
$$
D_{j,k}(uv)=
(D_{j+1,k+1}u)v+
u(D_{j+2,k+2}v).\eqno(1.32)$$
We first calculate
$$
g_{j}d_{k}(x,y)g_{j}^{-1}(uv)=
g_{j}d_{k}(x,y)\{(g_{j+1}^{-1}u)
(g_{j+2}^{-1}v)\}=
$$
$$
g_{j}(\{d_{k+1}(x,y)(g_{j+1}^{-1}u)\}
(g_{j+2}^{-1}v))+
g_{j}((g_{j+1}^{-1}u)
(d_{k+2}(x,y)g_{j+2}^{-1}v))
$$
$$=
g_{j+1}
(d_{k+1}(x,y)(g_{j+1}^{-1}u))g_{j+2}
(g_{j+2}^{-1}v))+
(g_{j+1}g_{j+1}^{-1}u)
(g_{j+2}d_{k+2}(x,y)g_{j+2}^{-1}v)$$
$$
=g_{j+1}(d_{k+1}(x,y)g_{j+1}^{-1}u)v+
u(g_{j+2}d_{k+2}(x,y)
g_{j+2}^{-1}v).$$
Similarly,
we find
$$
\{d_{k}(g_{j-k}x,g_{j-k}y)\}
(uv)=
$$
$$
\{d_{k+1}(g_{j-k}x,g_{j-k}y)u\}v+
u\{d_{k+2}(g_{j-k}x,g_{j-k}y)v\}
$$
and these prove the validity of Eq.(1.32),
since 
$$
g_{j}d_{k}(x,y)g_{j}^{-1}(uv) -d_{k}(g_{j-k}x,g_{j-k}y)(uv)=
$$
$$
(D_{j+1,k+1}u)v*u(D_{j+2.j*2}v).
$$
\par
If we choose $k=3$
in Eq.(1.32),
we obtain
$$
D_{j,3}(uv)=
(D_{j+1,1}u)v+
u(D_{j+2,2}v)=0$$
which gives 
$D_{j,3}=0,$
provided that the condition (B) holds.
\par
\noindent
On the other side,
the choice of $k=1$
or
$k=2$ in Eq.(1.32)
yields
$$
0=u(D_{j+2,3}v)=
(D_{j+1,3}u)v$$
for any $j=1,2,3$
and for any
$u,v\in A.$
Therefore,
it also gives
$D_{j,3}=0$
under the condition (C) .
This completes the proof.
$\square$
\par
We are now in position that we can construct a 
invariant sub-group of
$Trig(A)$ as follow:
\par
Let $A$ be a regular triality algebra satisfying 
either condition (B) or (C),
and set
$$
L_{0}=span <d_{j}(x,y),\forall x,y\in A,\forall j=1,2,3>.\eqno(1.33)$$
Then $L_{0}$ is a Lie algebra by Eq.(1.29b).
Moreover,
it is an ideal of the larger Lie algebra
$s\circ Lrt(A)$
by
Eq.(1.29a).
For any basis
$e_{1},e_{2} \cdots,e_{N}$
of $A$ with
$N=Dim\ A,$
and for any $\alpha _{j,\mu,\nu}\in F (j=1,2,3, \ \mu,\ \nu =1,2,\cdots, N),$
we set
$$
D_{j}=
\sum^{3}_{k=1}\sum^{N}_{\mu,\nu=1}\alpha_{j-k,\mu,\nu}
d_{k}(e_{\mu},e_{\nu})\in L_{0}\eqno(1.34)$$
for $j=1,2,3.$
Then,
$D=(D_{1},D_{2},D_{3})$
is a member of
$s\circ Lrt(A)$
by Eq.(1.29).
Thereofore its exponential
map
$$\xi_{j}=exp\ D_{j}\ (j=1,2,3)\eqno(1.35)$$
satisfies
$\xi_{j}(xy)=
(\xi_{j+1}x)(\xi_{j+2}y),$
i.e.,
$(\xi_{1},\xi_{2},\xi_{3})\in { Trig}(A),$
provided that
the exponential map is well-defined.
Moreover, for any $g=(g_{1},g_{2},g_{3})\in { Trig}(A).$
We calculate
$$
g_{j}D_{k}g_{j}^{-1}=
\sum^{3}_{l=1}\sum^{N}_{\mu,\nu=1}\alpha_{k=l,\mu,\nu}
\{d_{l}(g_{j-l}e_{\mu},e_{\nu})
+d_{j}(e_{\mu},g_{j-l}e_{\nu})\}
\in L_{0}\eqno(1.36)$$
by Eq.(1.30),
so that
$$
g_{j}\xi_{k}g_{j}^{-1}=
g_{j}(exp\ D_{k})g_{j}^{-1}=
exp(g_{j}D_{k}g_{j}^{-1})
\in exp\ L_{0}.\eqno(1.37)$$
Thereofre,
the group
$G_{0}$ generated by $<\xi_{1},\xi_{2},\xi_{3}>$
is a invariant sub-group of
$Trig (A)$.
\par
\vskip 3mm
{\bf Remark 1.4}
\par
\vskip 3mm
If a regular triality algebra satisfies Eqs.(1.27)
as well as
$$
d_{3}(x,y)z+
d_{3}(y,z)x+
d_{3}(z,x)y=0,\eqno(1.38)$$
then
$A$ is known as a pre-normal triality algebra([K-O.15]).
Moreover,
if we have
$$
Q(x,y,z)=d_{1}(z,xy)+
d_{2}(y,zx)+d_{3}(x,yz)=0\eqno(1.39)
$$
in addition,
$A$ is called a normal triality algebra in
[K-O.15] also.
\par
Further, suppose that $A$ is involutive with the involution map
$x\rightarrow{\bar x}$
as in Eq.(1.10)
and intoroduce the second bi-linear product in the same vector space of $A$ by
$$
x*y=\overline{ xy}=
{\bar y}{\bar x}.\eqno(1.40)$$
Then,
the resulting algebra
$A^{*}$
is also involutive, this algebra $A^{*}$
is called a conjugation algebra of $A$.
Furthermore, we note that a structurable algebra introduced by Allison([A.78],[A-F.93]) is
in a class of the conjugation algebra $A^{*}$ of the normal triality algebra $A$  and this nonassociative algebra $A^{*}$ contains 
associative, commutative Jordan and alternative algebras([K-O.15]).
\par
For a  Lie algebra $(A,\ [\ ,\ ])$, 
setting ${\bar x}=-x$, from $x*y=\overline{[x,y]}$,
we have
$$
\overline{[x,y]}=-[x,y]=[y,x]=[{\bar y},{\bar x}],$$
$$\overline{x*y}=[x,y]=\overline{[{\bar y}, {\bar x}]}={\bar y}*{\bar x}.
$$ 
Hence these  imply that 
$A$ and $A^{*}$ are involutive with respect to  the ${\bar x}=-x$.
\par 
If $A^{*}$
is unital (i.e., there is an element $e \in A^{*}$ such that $e*x=x*e=x$ for any $x\in A^{*}$)  and $A$ is a normal triality algebra,
then this algebra $A^{*}$ is structurable 
and both conditions (B) and (C)
are automatically satisfied (to see, [K-O.14]).
\par
Finally,
Eq.(1.1)
is rewritten as a modified global triality relation of
$$
{\bar g}_{j}(x*y)=
(g_{j+1}x)*(g_{j+2}y).\eqno(1.41)$$
Reviewing these results,
we can construct a invariant sub-group of ${ Trig} (A)$
from any regular triality algebra $A$ satisfying
condition (B) or
(C), in particular  from any structurable algebra $A^{*}$.
However,
for some cases,
we can more directly construct some invariant sub-group of $Trig (A)$, as it happenes for the case of $A$ being a 
symmetric composition algebra.
The case will be explored  in the next section,
while we will study its relation to the
corresponding
$s\circ Lrt(A)$
in section 3.
If we further restrict ourselves to the automorphism group,
then we can also construct some explict  automorphism of
 the Cayley algebra.
which we will study in section 4.
Finally,
other case of $Trig(A)$ for several algebras $A$ with involution 
 will be given in section 5.
\par

\vskip 5mm
{\bf 2\ Symmetric Composition Algebra}
\par
\vskip 3mm
This section is devoted for constructions of
${ Trig}(A)$ for symmetric composition algebras.
Let $A$ be an algebra over a field $F$ of characteristic not 
$2$ with a symmetric bi-linear non-degenerate
form $<\cdot|\cdot>.$
Suppose that it satisfies
$$
(xy)x=x(yx)=<x|x>y\eqno(2.1)$$
for any $x,y\in A.$
Then $A$ is called a symmetric composition algebra,
since it also satisfies the composition law
$$
<xy|xy>=<x|x><y|y>\ {\rm and}\ 
<xy|z>=<x|yz>.\eqno(2.2)$$
Conversely,
the validity of Eq.(2.2)
implies that of Eq.(2.1)
([O,95]).
\par
Linearizing Eq.(2.1),
it gives
$$
(xy)z+(zy)x=
x(yz)+z(yx)=
2<x|z>y\eqno(2.3)$$
for any $x,y,z\in A.$
Moreover if we replace
$z$
by $yz$
in the first relation of
Eq.(2.3) and note Eq.(2.2),
it yields
$$
(xy)(yz)=
2<x|yz>y-<y|y>zx.\eqno(2.4)$$
Also,
for any $g=(g_{1},g_{2},g_{3})\in { Trig}(A),$
we have
$$
<g_{j}x|g_{j}y>=<x|y>\ (j=1,2,3)\eqno(2.5)$$
by the following reason.
Applying
$g_{j}$
to both sides of
Eq.(2.1),
it gives
$$
<x|x>g_{j}y=
g_{j}\{x(yx)\}=
(g_{j+1}x)(g_{j+2}(yx))=$$
$$
(g_{j+1}x)\{(g_{j}y)(g_{j+1}x)\}=
<g_{j+1}x|g_{j+1}x>g_{j}y$$
so that we have 
$<x|x>=<g_{j+1}x|g_{j+1}x>.$
Linearizing the relation,
and letting
$j\rightarrow j-1,$
we obtain then Eq.(2.5).
Similarly,
any $t_{j}\in s\circ Lrt(A)$ satisfies
$$
<t_{j}x|y>=
-<x|t_{j}y>\ (j=1,2,3),\eqno(2.6)$$
by the relations  $t_{j}(x(yx))=<x|x>t_{j}(y)$ and  $<x|t_{j}(x)>=0$.
\par
Any symmetric composition algebra is known
(see [O-O.81],[E.97],
[K-M-R-T.98])
to be either a para-Hurwitz algebra or
$8$-dimensional pseudo-octonion algebra
where the para-Hurwitz algebra is the conjugate algebra of Hurwitz algebra with the para-unit $e$ satisfying
$ex=xe={\bar x}=2<e|x>e-x.$
That is, we emphasize that quaternion and Cayley algebras
with the product $x*y$ have a structure of a 
symmetric composition algebra with respect to new product $xy=\overline{x*y}$ , where ${\bar x}$ is the standard conjugation of $x$. 
Moreover,
the symmetric composition algebra is  a normal
triality algebra
([O.05]).
Further, 
conditions (B)and (C) of section 1 are 
automatically satisfied by this algebra.
Since $<\circ|\circ>$ is assumed to be non-degenerate,
there exists
$x\in A$
satisfying
$<x|x>\not= 0.$
Now Eq.(2.1)
implies
$<x|x>y\in AA$
for any element
$y\in A$,
so that we have
$AA=A.$
If $b\in A$
satisfies $bA=0.$
then Eq.(2.1)
also leads to
$<x|x>b=0$
and here
$b=0$.
Therefore both conditions (B) and (C) are
automatically satisfied.
Especially,
any symmetric composition
algebra is a regular
triality algebra satisfying both 
conditions (B) and (C)
(see Remark 1.4).
Thus,
we can construct some invariant sub-group
of ${ Trig}(A)$
as has been demonstrated in the previous sections.
However,
for special cases of 
$Dim\ A=1$
and $Dim\ A=2$,
we can construct entire
${ Trig}(A)$ as follows.
\par
\vskip 3mm
{\bf Example 2.1.\ Case of $Dim\ A=1$}
\par
\vskip 3mm
We write
$A=Fe$,
where $e\in A$
satisfies
$ee=e$
with $<e|e>=1.$
Then for any
$g=(g_{1},g_{2},g_{3})\in {\rm Trig}(A),$
we can write
$$
g_{j}(e)=\alpha_{j}e,\ (j=1,2,3)$$
for some
$\alpha_{j}=\alpha_{j+1}\alpha_{j+2}$
while Eq.(2.5) yields
$(\alpha_{j})^{2}=1.$
Therefore,
there gives
$$
\alpha_{1}^{2}=\alpha_{2}^{2}=\alpha_{3}^{2}=1\ 
{\rm and}\ 
\alpha_{1}\alpha_{2}\alpha_{3}=1$$
so that 
${\rm Trig}(Fe)$ is isomorphic to the Klein's
$4$-group $K_{4}$.
Note that we have
${\rm Auto}(Fe)=<Id>$
being trivial.
\par
\vskip 3mm
{\bf Example 2.2.\ Case of $Dim\ A=2$}
\par
\vskip 3mm
In this case,
$A={\rm span}<e,f>$
where
$e, f\in A$
satisfy
$$
ee=e,\ ff=-e,\ ef=fe=-f\eqno(2.7a)$$
$$
<e|e>=<f|f>=1,\ <e|f>=<f|e>=0.\eqno(2.7b)$$
Note that $A$ is Abelian,
and  any
$x,y\in A$
satisfies the quadratic
relation of
$$
xy=yx=-<e|x>y-
<e|y>x+
(4<e|x><e|y>-
<x|y>)e.\eqno(2.8)$$
Note also that $e$ is the para-unit of $A$,
i.e.,
$$
ex=xe={\bar x}=
2<e|x>e-x.\eqno(2.9)$$
\par
Let us write
$$
p_{j}=g_{j}(e),q_{j}=g_{j}(f),\ (j=1,2,3).\eqno(2.10)$$
Then,
the global triality relation Eq.(1.1)
is equivalent to the validity of
$$
p_{j+1}p_{j+2}=p_{j},\ p_{j+1}q_{j+2}
=q_{j+1}p_{j+2}=-q_{j},\ q_{j+1}q_{j+2}=-p_{j}\eqno(2.11)$$
for $j=1,2,3$
where the indices over $j$ are defined modulo $3$
i.e.,
$p_{j\pm3}=p_{j},q_{j\pm3}=q_{j}.$
Also,
Eqs.(2.5) and (2.7b) give
$$
<p_{j}|p_{j}>=
<q_{j}|q_{j}>=1,\ <p_{j}|q_{j}>=0,\ (j=1,2,3).\eqno(2.12)$$
\par
The general solution of Eqs.(2.11) and
(2.12) is given by the following Proposition.
\vskip 3mm
{\bf Proposition 2.3}
\it
\par
\vskip 3mm
Let $q_{j}\in A$ for
$j=1,2,3$
be any three
element of $A$
satisfying
$$
<q_{1}|q_{1}>=<q_{2}|q_{2}>=
<q_{3}|q_{3}>=1,\ {\rm and}\ 
<q_{3}|q_{1}q_{2}>=0.\eqno(2.13)$$
If we introduce $p_{j}$
by 
$$
p_{j}=-q_{j+1}q_{j+2},\eqno(2.14)$$
then they satisfy Eqs.(2.11) and (2.12).
Conversly,
if $p_{j}$
and
$q_{j}$
satisfy Eqs>(2.11)
and
(2.12),
then
Eq.(2.13)
holds.
\par
\vskip 3mm
\rm
{\bf Proof}
\par
\vskip 3mm
We first note that
$<q_{3}|q_{1}q_{2}>=0$
implies the validity also of
$$
<q_{j}|q_{j+1}q_{j+2}>=0\eqno(2.15)$$
for all
$j=1,2,3,$
because of the associate law
$<xy|z>=<x|yz>.$
\par
Suppose now that Eq.(2.13)
with
Eq.(2.14)
holds:
We calculate then
$$
p_{j+1}q_{j+2}=
-(q_{j+2}q_{j})q_{j+2}=
-<q_{j+2}|q_{j+2}>q_{j}=-q_{j}$$
$$
q_{j+1}p_{j+2}=
-q_{j+1}(q_{j}q_{j+1})=
-<q_{j+1}|q_{j+1}>q_{j}=-q_{j}.$$
Finally,
we compute
$$
p_{j+1}p_{j+2}=
(q_{j+2}q_{j})
(q_{j}q_{j+1})
=
2<q_{j+2}|q_{j}q_{j+1}>
q_{j}-<q_{j}|q_{j}>q_{j+1}q_{j+2}$$
by Eq.(2.4)
so that we find
$$
p_{j+1}p_{j+2}=
-q_{j+1}q_{j+2}=p_{j}$$
because of Eq.(2.14),
proving
Eq.(2.11).
Moreover
$$
<p_{j}|p_{j}>=
<q_{j+1}q_{j+2}|q_{j+1}q_{j+2}>=
<q_{j+1}|q_{j+1}>
<q_{j+2}|q_{j+2}>=1$$
and
$$
<q_{j}|p_{j}>=
-<q_{j}|q_{j+1}q_{j+2}>=0$$
which establishes Eq.(2.12).
Conversely,
if
$<p_{j}|q_{j}>=0$
with
$p_{j}=-q_{j+1}q_{j+2},$
then they lead to
$<q_{j}|q_{j+1}q_{j+2}>=0$
or
$<q_{3}|q_{1}q_{2}>=0.\square$
\par
\vskip 3mm
{\bf Remark 2.4(a)}
\par
\vskip 3mm
If we write
$$q_{j}=\mu_{j}e+\nu_{j}f,\ (j=1,2,3)
\eqno (2.16a)
$$
for
$\mu_{j},\nu_{j}\in F\ (j=1,2,3),$
then
Eq.(2.13)
gives
$$
\mu_{j}^{2}+\nu_{j}^{2}=1,\ 
\mu_{1}\mu_{2}\mu_{3}=\mu_{1}\nu_{2}\nu_{3}+
\mu_{2}\nu_{3}\nu_{1}+
\mu_{3}\nu_{1}\nu_{2}\eqno(2.16b)$$
while Eq. (2.14) for
$$
p_{j}=\alpha_{j}e+\beta_{j}f\eqno(2.17)$$
these will leads to
$$
-\alpha_{j}=\mu_{j+1}\mu_{j+2}-\nu_{j+1}\nu_{j+2},\eqno(2.18a)$$
$$
-\beta_{j}=-(\mu_{j+1}\nu_{j+2}+
\nu_{j+1}\mu_{j+2}).\eqno(2.18b)$$
Especially,
Eq.(2.16)
allows two indepent variables
in whose terms
other
$\mu_{j}$
and
$\nu_{j}$
are determined.
Thus,
if the field $F$ is either real or
complete filed,
then the resultained group manifold
${ Trig}A$
will become some
$2$-dimensional real or complex manifold.
However,
when
$F$ is a finite
field,
then
${ Trig}(A)$ will become to a finite goup.
Also,
it may be of some
interest to see that
$\alpha_{j}=<e|p_{j}>$
and
$\beta_{j}=<f|p_{j}>$
satisfy
polynomial identities of
$$
2\alpha_{1}\alpha_{2}\alpha_{3}-
(\alpha_{1}^{2}+\alpha_{2}^{2}+\alpha_{3}^{2})+
1=0,\eqno(2.19a)$$
$$
(\beta_{1}^{2}+\beta_{1}^{2}+\beta_{3}^{2})^{2}+
4(\beta_{1}\beta_{2}\beta_{3})^{2}=
4\{\beta_{1}^{2}\beta_{1}^{2}+
\beta_{2}^{2}\beta_{3}^{2}+
\beta_{3}^{2}\beta_{1}^{2}\}.\eqno(2.19b)$$
However,
since its relevance is not clear,
we will not go into its detail.
\par
Next we consider
for $s\circ Lrt(A),$
any
$t=(t_{1},t_{2},t_{3})\in s\circ Lrt(A)$
is given similarly by
$$
t_{j}e=\lambda_{j}f,\ t_{j}f=-\lambda_{j}e\ (j=1,2,3)\eqno(2.20a)$$
for some
$\lambda_{j}\in F$
satisfying
$$
\lambda_{1}+\lambda_{2}+\lambda_{3}=0.\eqno(2.20b)
$$
\par
\vskip 3mm
{\bf Remark 2.4 (b)}\ $  Auto\ (A)$ with $Dim \ A=2$ 
\par
\vskip 3mm
In contrast to ${ Trig}(A),\ 
{ Auto}(A)$
for the case of
$Dim\ A=2$
turns out to be always
a finite group of either
$Z_{2}$
if
$\sqrt{3}\notin F$
or
$S_{3}$
if 
$\sqrt{3}\in F.$
Note that
${ Auto}(A)$
is a special case of
${ Trig}(A)$
for
$g_{1}=g_{2}=g_{3}(=\sigma)$
with
$$
\sigma(xy)=(\sigma x)(\sigma y),\eqno(2.21a)$$
Setting $p_{1}=p_{2}=p_{3}=p,$ and $q_{1}=q_{2}=
q_{3}=q,$ then we have
$$
\sigma(e)=p,\ \sigma(f)=q\eqno(2.21b)
$$
and Eq.(2.13) becomes
$$
<q|q>=1,\ <q|q^{2}>=0\eqno(2.22)$$
while Eq.(2.14)
determines
$p$ by
$$
p=-qq.\eqno(2.23)$$
Setting $x=y=q$
in Eq.(2.8),
we find
$$
q^{2}=-2<e|q>q+
(4<e|q><e|q>-
<q|q>)e\eqno(2.24)$$
so that
Eq.(2.22)
gives
$$
0=<q|q^{2}>=
-2<e|q><q|q>+
(4<e|q><e|q>-
<q|q>)<e|q>$$
$$
=<e|q>\{4<e|q><e|q>-3\},\eqno(2.25)$$
whose solutions are
\par
(i)
$$
<e|q>=0\eqno(2.26a)$$
or
\par
(ii)
$$
<e|q>=\pm
{{\sqrt{3}}\over 2},\eqno(2.26b)$$
assuming
$\sqrt{3}\in F$
for the
2nd solution if
Eq.(2.26b).
\par
In terms of $e$ and $f$,
then these imply.
\par
\vskip 3mm
{\bf Case 1}
\par
\quad
$\sigma:e\rightarrow e,\ f\rightarrow \epsilon f,\ (\epsilon=\pm 1)$
\par
\vskip 3mm
{\bf Case 2 (if $\sqrt{3}\in F$)}
\par
\quad
$
\sigma:e\rightarrow {1\over 2}(-e+\sqrt{3}\epsilon_{1}\epsilon_{2}f_{1})$
\par
\quad\quad
$f\rightarrow {1\over 2}(
\sqrt{3}\epsilon_{1}e+\epsilon_{2}f),$
\par
\quad\quad where
$\epsilon_{1}^{2}=\epsilon_{2}^{2}=1.$
\par
\vskip 3mm
Let $P$ and
$Q\in {\rm Auto}(A)$
be given by
\par
\quad 
$P:e\rightarrow e,\ f\rightarrow -f$
\par
\quad
$Q:e\rightarrow
{1\over 2}(-e+\sqrt{3}f)$,  
\quad
$f\rightarrow {1\over 2}(-\sqrt{3}e-f).$
\par
Then, they satisfy
$$
P^{2}=Q^{3}=Id,\ QPQ=P$$
which generate the symmetric group $S_{3}$
(or equivalently the
dihedral group $D_{3}$).
Thus, we conclude
$$
{ Auto}(A)=
\left\{
\begin{array}{l}
Z_{2}=<Id,P>,\ {\rm if}\sqrt{3}
\not\in F\\
S_{3}=<Id,P,Q,Q^{2},PQ,PQ^{2}>,\ {\rm if}\ \sqrt{3}\in F,
\end{array}
\right.
$$
which are finite group in contrast to the case
of
${ Trig}(A)$
(at least for real or
complex field $F$).
Correspondly,
we have
${ Der}(A)=0,$
for the derivation assuming $3\not=0.$
\par
Other cases of
${ Dim}\  A=4$
(para-quotanion algebra) and
${ Dim} \ A=8$
(either para-or pseudo octonion algebra)
are unfortunately,
too complicated to be discussed in a similar way.
However,
we can still construct some
sub-groups of
${ Trig}(A)$ as follows.
For any two element $a,b\in A$
satisfying
$<a|a>=<b|b>=1,$
we set
$a_{1}=a,\ a_{2}=b$
and define
by $a_{3}=a_{1}a_{2}=ab.$
Then,
we find
$$
a_{j}a_{j+1}=a_{j+2},\ <a_{j}|a_{j}>=1\eqno(2.27)$$
for any
$j=1,2,3,$
where
the indices over $j$ are defined modulo
$3$,i.e.,
$$
a_{j\pm3}=a_{j}.$$
For example,
we calculate
$$
a_{3}a_{1}=
(a_{1}a_{2})a_{1}=
<a_{1}|a_{1}>a_{2}=a_{2},$$
$$
<a_{3}|a_{3}>=
<a_{1}a_{2}|a_{1}a_{2}>=
<a_{1}|a_{1}><a_{1}|a_{2}>=1$$
etc.
\par
In view of this, we introduce
$$
\Sigma=\{a=(a_{1},a_{2},a_{3})\in A^{3}|
a_{j}a_{j+1}=a_{j+2},\ 
<a_{j}|a_{j}>=
1,\ 
\forall j=1,2,3\}.\eqno(2.28)
$$
We now prove
\par
\vskip 3mm
{\bf Theorem 2.5}
\par
\it
\vskip 3mm
Let $A$ be a symmetric composition algebra.
For any $a=(a_{1},a_{2},a_{3})\in \Sigma,$
we set
$$
\sigma_{j}(a)=
R(a_{j+1})R(a_{j+2}),\eqno(2.29a)$$
$$
\theta_{j}(a)=L(a_{j+2})L(a_{j+1})\eqno(2.29b)$$
for
$j=1,2,3.$
\par
Then,
they satisfy
\par
(i)
$$
\sigma_{j}(a)(xy)=
(\sigma_{j+1}(a)x)(\sigma_{j+2}(a)y)\eqno(2.30a)$$
\par
(ii)
$$\theta_{j}(a)(xy)=
(\theta_{j+1}(a)x)(\theta_{j+2}(a)y)\eqno(2.30b)$$
\par
(iii)
$$
\sigma_{j}(a)\theta_{j}(a)=
\theta_{j}(a)\sigma_{j}(a)=Id\eqno(2.30c)$$
\par
(iv)
$$
\theta_{j}(a)\theta_{j+1}(a)\theta_{j+2}(a)=
\sigma_{j+2}(a)\sigma_{j+1}(a)\sigma_{j}(a)=Id\eqno(2.30d)$$
\par
(v)
$$
<\sigma_{j}(a)x|y>=
<x|\theta_{j}(a)y>\eqno(2.30e)$$
\par
(vi)
$$
<\sigma_{j}(a)x|
\sigma_{j}(a)y>=
<\theta_{j}(a)x|
\theta_{j}(a)y)>=
<x|y>.\eqno(2.30f)$$
Especially,
both triples
$\sigma(a)=
(\sigma_{1}(a),\sigma_{2}(a),\sigma_{3}(a))$
and
$\theta(a)=
(\theta_{1}(a),\theta_{2}(a),\theta_{3}(a))$
are elements of ${ Trig}(A).$
\par
\vskip 3mm
\rm
{\bf Proof}
\par
\vskip 3mm
We will first prove
(iii)-(vi)
as follows:
Eq.(2.1)
is rewritten as
$$
R(x)L(x)=
L(x)R(x)=<x|x>Id\eqno(2.31)$$
so that we calculate
$$
\sigma_{j}(a)\theta_{j}(a)~
=R(a_{j+1})R(a_{j+2})
L(a_{j+2})L(a_{j+1})=
R(a_{j+2})L(a_{j+1})=Id,$$
$$
\theta_{j}(a)\sigma_{j}(a)=
L(a_{j+2})L(a_{j+1})
R(a_{j+1})R(a_{j+2}=
L(a_{j+2})R(a_{j+2})=Id$$
which prove
Eq.(2.30c).
We next also have 
$$
<\sigma_{j}(a)x|y>=
<R(a_{j+1})R(a_{j+2})x|y>
=<(xa_{j+2})a_{j+1}|y>$$
$$
=<xa_{j+2}|a_{j+1}y>=
<x|a_{j+2}(a_{j+1}y)>=
<x|L(a_{j+2}L(a_{j+1})y)=
<x|\theta_{j}(a)y>$$
which is Eq.(2.30e).
Then,
setting
$y\rightarrow \theta_{j}(a)y$
in
Eq.(2.30e)
and noting
Eq.(2.30c),
we obtain Eq.(2.30f).
\par
We next rewrite Eqs.(2.29)
as
$$
\sigma_{j}(a)(x)=
2<a_{j+1}|x>a_{j+2}-a_{j}x\eqno(2.32a)$$
$$
\theta_{j}(a)x=
2<a_{j+2}|x>a_{j+1}-xa_{j}\eqno(2.32b)$$
since we have
$$
\sigma_{j}(a)x=
(xa_{j+2})a_{j+1}=
2<x|a_{j+1}>
a_{j+2}-
(a_{j+1}a_{j+2})x=
2<a_{j+1}|x>a_{j+2}-a_{j}x$$
and
$$
\theta_{j}(a)x=
a_{j+2}(a_{j+1}x)=
2<a_{j+2}|x>a_{j+1}-
x(a_{j+1}a_{j+2})=
2<a_{j+2}|x>a_{j+1}-xa_{j}$$
in view of
Eq.(2.3) and
(2.28).
We then calculate
$$
\sigma_{j+2}(a)\sigma_{j+1}(a)x=
2<a_{j}|\sigma_{j+1}(a)x>a_{j+1}-
a_{j+2}(\sigma_{j+1}(a)x)$$
$$
=2<a_{j}|2<a_{j+2}|x>a_{j}-a_{j+1}x>a_{j+1}
-a_{j+2}\{2<a_{j+2}|x>a_{j}-a_{j+1}x\}$$
$$
=\{
4<a_{j+2}|x>
<a_{j}|a_{j}>
-2<a_{j}a_{j+1}|x>\}
a_{j+1}
-2<a_{j+2}|x>a_{j+2}a_{j}
+a_{j+2}(a_{j+1}x)
$$
$$
=\{4<a_{j+2}|x>-
2<a_{j+2}|x>\}
a_{j+1}
-2<a_{j+2}|x>
a_{j+1}+
a_{j+2}(a_{j+1}x)=
a_{j+2}(a_{j+1}x)
$$
$$=
\theta_{j}(a)x,
$$
so that
$\sigma_{j+2}(a)\sigma_{j+1}(a)=
\theta_{j}(a).$
Multiplying
$\sigma_{j}(a)$ from the left and noting
Eq.(2.30c),
we obtain
$\sigma_{j+2}(a)\sigma_{j+1}(a)\sigma_{j}(a)= Id.$
If we further multiply
$\theta_{j}(a)\theta_{j+1}(a)\theta_{j+2}(a)$
from the left to this,
it also reproduce
$\theta_{j}(a)\theta_{j+1}(a)\theta_{j+2}(a)=Id,$
proving Eqs.(.2.30c).
\par
It now remains to prove Eqs.(2.30a)
and (2.30b).
We calculate
$$
(\sigma_{j}(a)x)
(\sigma_{j+1}(a)y)=
(\sigma_{j}(a)x)
\{(ya_{j})a_{j+2}\}
=2<\sigma_{j}(a)x|a_{j+2}>ya_{j}
-a_{j+2}\{(ya_{j})(\sigma_{j}(a)x)\}$$
$$
=2<x|\theta_{j}(a)a_{j+2}>ya_{j}
-a_{j+2}\{(ya_{j})(2<x|a_{j+1}>a_{j+2}
-a_{j}x)\}$$
$$
=2<x|a_{j+2}(a_{j+1}a_{j+2})>ya_{j}
-2<x|a_{j+1}>a_{j+2}
\{(ya_{j})a_{j+2}\}
+a_{j+2}\{(ya_{j})(a_{j}x)\}$$
$$
=2<x|<a_{j+2}|a_{j+2}>a_{j+1}>ya_{j}-2<x|a_{j+1}>ya_{j}
+a_{j+2}\{2<y|a_{j}x>a_{j}-<a_{j}|a_{j}>xy\}
$$
$$
=2<xy|a_{j}>
a_{j+1}-
a_{j+2}(xy)=
\sigma_{j+2}(a)(xy),
$$
where we utilized Eq.(2.4).
This proves Eq.(2.30a).
If we replace
$x$ and $y$ by
$\theta_{j+1}(a)x$
and $\theta_{j+1}(a)y$
respectively
in Eq.(2.30a),
we also find Eq.(2.30b).
This completes the proof. $\square$
\par
Noting
$$
\sigma(a)=(\sigma_{1}(a),\sigma_{2}(a),\sigma_{3}(a))
\in { Trig}(A)\eqno(2.33a)$$
$$
\theta(a)=
(\theta_{1}(a),
\theta_{2}(a),
\theta_{3}(a))
\in
{ Trig}(A)\eqno(2.33b)$$
by Theorem 2.5,
these generate a sub-group of
${ Trig}(A)$.
However except for the cases of
${ Dim}A=1$
and
${ Dim}A=2,$
it is not a invariant sub-group of
${ Trig}(A)$
by the following reason.
For any
$g=(g_{1},g_{2},g_{3})
\in{\rm Trig}(A),$
Eqs.(1.25)
yield
$$
g_{j}L(x)L(y)=
L(g_{j+1}x)g_{j+2}L(y)=
L(g_{j+1}x)L(g_{j}y)g_{j+1},\eqno(2.34a)$$
$$
g_{j}R(x)R(y)=
R(g_{j+2}x)g_{j+1}R(y)=
R(g_{j+2}x)R(g_{j}y)g_{j+2}\eqno(2.34b)$$
so that
$$
g_{j}\theta_{k}(a)=
L(g_{j+1}a_{k+2})
L(g_{j}a_{k+1})g_{j+1}\eqno(2.34c)$$
$$
g_{j}\sigma_{k}(a)=
R(g_{j+2}a_{k+1})
R(g_{j}a_{k+2})g_{j+2}.\eqno(2.34d)$$
\par
We next note
$$
(g_{j}a_{j})(g_{j+1}a_{j+1})=
g_{j+2}(a_{j}a_{j+1})=
g_{j+2}a_{j+2}$$
so that we have
$$
ga=(g_{1}a_{1},g_{2}a_{2},g_{3}a_{3})\in \Sigma
\eqno(2.35)
$$
Moreover,
we have similarly
$$
\phi g=
	(g_{2},g_{3},g_{1})\in { Trig}(A)
$$
$$
\phi^{2}g=
(g_{3},g_{1},g_{2})\in
{ Trig}(A)
$$
for
$\phi:1\rightarrow 2\rightarrow 3\rightarrow 1.$
Then,
Eqs.(2.34)
are rewritten as
$$
g_{j}\theta_{k}(a)g_{j+1}^{-1}=
\theta_{k}((\phi^{(j-k-1)}g)a)\eqno(2.36a)$$
$$
g_{j}\sigma_{k}(a)
g_{j+2}^{-1}=
\sigma_{k}((\phi^{(j-k+1)}g)a)\eqno(2.36b)$$
for any
$j,k=1,2,3.$
Especially,
for
$j=k,$
we obtain
$$
g_{j}\theta_{j}(a)g_{j+1}^{-1}=
\theta_{j}((\phi g)a)\eqno(2.37a)$$
$$
g_{j}\sigma_{j}(a)g_{j+2}^{-1}=
\sigma_{j}((\phi^{2}g)a).\eqno(2.37b)
$$
\par
Since the left sides of Eqs.(2.37) are not of a form
$g_{j}\times g_{j},$
the
sub-group $G_{0}$
generated by
$\sigma_{j}(a)$
and $\theta_{j}(a)$
is not in general invariant under 
${ Trig}(A).$
However, cases of
${ Dim}A=1$
and
${ Dim}\ A=2$ are exceptional,
since they are Abelian
and hence satisfy
$R(x)=L(x).$
Then,
we can rewrite
$$
\sigma_{j}(a)=R(a_{j+1})R(a_{j+2})=
R(a_{j+1})L(a_{j+2})$$
$$
\theta_{j}(a)=
L(a_{j+2})L(a_{j+1})=
L(a_{j+2})R(a_{j+1}).$$
Then,
Eqs.(1.26)
show that these are invaiant under
${ Trig}(A).$
\par
Returning to the original problem,
we can construct some invariant
sub-group of
${ Trig}(A).$
Let
$G$ be a sub-group generated now by
$$\sigma_{j}(a)\theta_{j}(b),\ \ 
and\ \ 
\theta_{j}(a)\sigma_{j}(b),\eqno(2.38a)$$
or
$$\sigma_{j}(a)\sigma_{j}(b)\sigma_{j}(c),\ \ 
and\ \ 
\theta_{j}(a)\theta_{j}(b)\theta_{j}(c)\eqno(2.38b)$$
for any
$a,b,c\in \Sigma.$
Then,
by Eq.(2.37),
$X=(X_{1},X_{2},X_{3})\in G$
will be transformed by
${ Trig}(A)$
into
$$
g_{j}X_{j}g_{j}^{-1}=
X_{j}^{'}$$
for another
$X^{'}=
(X_{1}^{'},X_{2}^{'},X_{3}^{'})\in G$
so that
$G$ is now a invariant sub-group of
${ Trig}(A).$

\par
In the final of this section,
for the 8 dimensional pseudo octonion algebra,  we note that
there is an example of elements in $\Sigma$ as follows.
\par
By means of the notations and  Eqs (4.26) in ([O.95]),
we have basis  $e_{1}, \cdots , e_{8}$  such that

$$
e_{j} e_{k}={\sum_{l=1}^8 ({\sqrt 3} d_{jkl} \mp f_{jkl})e_{l} }
\ \ and \ <e_{i}|e_{j}>= \delta_{ij}.
$$

For these special basis $ e_{1}, e_{2}, e_{3}$  in particular  the case of + sign in $\mp$, 
we may obtain the validity;
$$
e_{1} e_{2}=e_{3},\ \ e_{2} e_{3}=e_{1}, \ \ e_{3} e_{1}=e_{2},
$$
since $d_{12l}=d_{23l}=d_{31l}=0$, ($l=1,\cdots , 8$)
and $f_{123}=f_{231}=f_{312}=1$, $f_{12l} =0  \ (if \ l\ne 3)$,
$f_{23l}=0\ \ (if \ l\ne 1)$, $f_{ 31l }=0\ \ (if\ l \ne 2)$,
(for the numerical values of $d_{jkl} \ and\ f_{jkl}$, and 
an application to physics,
 see [G.62]). 

\par
Hence this implies that there is an element $(e_{1},e_{2},e_{3}) \in \Sigma $
satisfying the assumption in Theorem 2.5.
Furthemore, as this algebra is defined by 3x3 matrices, it is easy to see that
the  transpose of matrix induces an involution with respect to $\overline{x}=\ ^{t} x$.
 \vskip 5mm
{\bf 3\ Local Triality Relation}
\par
\vskip 3mm
In this section,
we will study the relationship between the invariant sub-group of 
${ Trig}(A)$ given in section 2 and 
its local corespondent
$s\circ Lrt(A)$ defined as in (1.18).
\par
From now on, we assume that the  algebra $A$  
is  a symmetric composition algebra over the field $F$
characteristic not 2 in Section 3 and Section 4
unless otherwise specifield.
\par 
Let
$$
\Sigma=
\{a=(a_{1},a_{2},a_{3})\in A^{3}|
a_{j}a_{j+1}=
a_{j+2},
<a_{j}|a_{j}>=1,\ \forall j=1,2,3\}\eqno(3.1)$$
as in Eq.(2.28).
For a given
$a=(a_{1},a_{2},a_{3})\in \Sigma,$
we introduce another triple
$\Lambda(a)$
by 
$$
\Lambda(a)=\{p=(p_{1},p_{2},p_{3})
\in A^{3}| a_{j}p_{j+1}+p_{j}a_{j+1}=p_{j+2}, <p_{j}|a_{j}>=0,\forall j=1,2,3\}\eqno(3.2)
$$
where the indices over
$j$ are defined modulo $3$
with
$p_{j+3}=p_{j}.$
We note then
that
$\Lambda(a)=
{\rm span}<p_{1},p_{2}>$
for
$p_{1},p_{2}$
satisfying
$<p_{j}|a_{j}>=0,$
only for
$j=1,2,$
as we see below.
For any such
$p_{1}$
and
$p_{2},$
we define
$p_{3}$ by
$$
p_{3}=a_{1}p_{2}+p_{1}a_{2},$$
we then calculate
$$
<a_{3}|p_{3}>=
<a_{3}|a_{1}p_{2}>+
<a_{3}|p_{1}a_{2}>=
<a_{3}a_{1}|p_{2}>+
<a_{2}a_{3}|p_{1}>
$$
$$
=<a_{2}|p_{2}>+
<a_{1}|p_{1}>=0.$$
$$
a_{2}p_{3}+p_{2}a_{3}=
a_{2}(a_{1}p_{2}+p_{1}a_{2})+
p_{2}(a_{1}a_{2})=
2<a_{2}|p_{2}>
a_{1}+<a_{2}|a_{2}>p_{1}=p_{1},$$
and
$$
a_{3}p_{1}+p_{3}a_{1}=
(a_{1}a_{2})p_{1}+
(a_{1}p_{2}+p_{1}a_{2})a_{1}$$
$$
=(a_{1}a_{2})p_{1}+
(p_{1}a_{2})a_{1}+
(a_{1}p_{2})a_{1}=
2<a_{1}|p_{1}>a_{2}+
<a_{1}|a_{1}>p_{2}=p_{2},$$
so that
$(p_{1},p_{2},p_{3})\in A(a).$
Note that
$\Lambda(a)$
is then a vector space
over the field
$F$
in contrast to
$\sum.$
\par
Moreover,
we further define
$q_{j}\in A$
by
$$
q_{j}=a_{j+1}p_{j+2}=p_{j}-p_{j+1}a_{j+2}.\eqno(3.3)$$
We then find
\par
\vskip 3mm
{\bf Lemma 3.1}
\par
\vskip 3mm
\it
Under the assumption
as in above, we have
$q=(q_{1},q_{2},q_{3})\in \Lambda(p).$
Moreover,
$p_{j}$
is computed conversely
by
$$
p_{j}=
q_{j+1}a_{j+2}=
q_{j}-a_{j+1}q_{j+2}.\eqno(3.4)
$$
\par
\vskip 3mm
\rm
{\bf Proof}
\par
\vskip 3mm
We first note
$$
<q_{j}|a_{j}>=
<a_{j+1}p_{j+2}|
a_{j+1}a_{j+2}>=
<a_{j+1}|a_{j+1}>
<p_{j+2}|a_{j+2}>=0.$$
Also,
we calculate
$$
q_{j+1}a_{j+2}=
(a_{j+2}p_{j+3})a_{j+2}=
<a_{j+2}|a_{j+2}>
p_{j+3}=p_{j}.$$
Finally,
we find
$$
a_{j}q_{j+1}+q_{j}a_{j+1}=
a_{j}(a_{j+2}p_{j+3})+
(a_{j+1}p_{j+2})a_{j+1}$$
$$
=2<a_{j}|p_{j}>a_{j+2}-
p_{j}(a_{j+2}a_{j})
+<a_{j+1}|a_{j+1}>p_{j+2}
$$
$$
=-p_{j}a_{j+1}+p_{j+2}=q_{j+2}.$$
These prove Lemmea 3.1.$\square$
\par
\vskip 3mm
{\bf Theorem 3.2}
\it
\par 
\vskip 3mm
For any
$a\in\Sigma$ and
$p\in \Lambda(a),$
if we introduce
$D_{j}(a,p)\in
{ End}\ A$
by
$$
D_{j}(a,p)x=
(p_{j+1}x)a_{j+1}+a_{j}(xq_{j}).\eqno(3.5)
$$
Then,
they satisfy
$$
D_{j}(a,p)(xy)=
(D_{j+1}(a,p)x)y+
x(D_{j+2}(a,p)y)\eqno(3.6)$$
so that
$$
D(a,p)=
(D_{1}(a,p),D_{2}(a,p),
D_{3}(a,p))
\in
s\circ Lrt(A).$$
\par
\rm
For a proof of this theorem,
we require
\par
\vskip 3mm
{\bf Lemma 3.3}
\par
\vskip 3mm
\it
Under the assumption as in Theorem 3.2, 
we have
$$
D_{j}(a,p)x=
2<q_{j+2}|x>a_{j+2}-
2<a_{j+2}|x>q_{j+2}+
a_{j}(xp_{j})\eqno(3.7a)$$
$$
=2<p_{j+2}|x>a_{j+2}-
2<a_{j+2}|x>p_{j+2}+
(q_{j+1}x)a_{j+1}.\eqno(3.7b)
$$
\par
\vskip 3mm
\rm
{\bf Proof}
\par
\vskip 3mm
We calculate
$$
a_{j}(xq_{j})=
a_{j}\{xp_{j}-x(p_{j+1}a_{j+2})\}
=a_{j}(xp_{j})-a_{j}\{2<x|a_{j+2}>p_{j+1}
-a_{j+2}(p_{j+1}x)\}
$$
$$
=a_{j}(xp_{j})-
2<x|a_{j+2}>a_{j}p_{j+1}+
a_{j}\{a_{j+2}(p_{j+1}x)\}$$
$$
=a_{j}(xp_{j})-
2<x|a_{j+2}>a_{j}p_{j+1}+
2<a_{j}|p_{j+1}x>a_{j+2}-
(p_{j+1}x)(a_{j+2}a_{j})$$
$$
=a_{j}(xp_{j})-
2<x|a_{j+2}>q_{j+2}+
2<a_{j}p_{j+1}|x>a_{j+2}-(p_{j+1}x)a_{j+1}$$
$$
=a_{j}(xp_{j})-2<x|a_{j+2}>q_{j+2}+
2<q_{j+2}|x>a_{j+2}
-(p_{j+1}x)a_{j+1}$$
so that we obtain
$$
D_{j}(a,p)x=a_{j}(xq_{j})+(p_{j+1}x)a_{j+1}$$
$$
=a_{j}(xp_{j})-2<x|a_{j+2}>q_{j+2}+
2<q_{j+2}|x>a_{j+2}$$
which is Eq.(3.7a).
In order to show
Eq.(3.7b),
we note
$$
a_{j}(xp_{j})=
a_{j}\{x(q_{j+1}a_{j+2})\}=
a_{j}\{2<x|a_{j+2}>q_{j+1}-
a_{j+2}(q_{j+1}x)\}$$
$$
=2<x|a_{j+2}>a_{j}q_{j+1}-
a_{j}\{a_{j+2}(q_{j+1}x)\}.\eqno(3.8)
$$
However, by $q_{j+1}=a_{j+2}p_{j}, $ 
we further compute
$$
a_{j}q_{j+1}=
a_{j}(a_{j+2}p_{j})=
2<a_{j}|p_{j}>
a_{j+2}-p_{j}(a_{j+2}a_{j})=
-p_{j}a_{j+1}$$
and
$$
a_{j}\{a_{j+2}(q_{j+1}x)\}=
2<a_{j}|q_{j+1}x>a_{j+2}-
(q_{j+1}x)(a_{j+2}a_{j})$$
$$
=2<a_{j}q_{j+1}|x>a_{j+2}-(q_{j+1}x)a_{j+1}$$
$$
=-2<p_{j}a_{j+1}|x>a_{j+2}-
(q_{j+1}x)a_{j+1}$$
so that
Eq.(3.8)
becomes
$$
D_{j}(a,p)x=
2<q_{j+2}|x>a_{j+2}-
2<a_{j+2}|x>q_{j+2}$$
$$
-2<x|a_{j+2}>p_{j}a_{j+1}+
2<p_{j}a_{j+1}|x>
a_{j+2}+(q_{j+1}x)a_{j+1}$$
$$
=2<(q_{j+2}+p_{j}a_{j+1})|x>a_{j+2}-
2<a_{j+2}|x>(q_{j+2}+p_{j}a_{j+1})+
(q_{j+1}x)a_{j+1}$$
$$
=2<p_{j+2}|x>a_{j+2}-
2<a_{j+2}|x>p_{j+2}+
(q_{j+1}x)a_{j+1}$$
which is Eq.(3.7b).
This completes the proof of Lemma 3.3$\square$
\par
We are now in a position to prove Theorem 3.2.
We calculate
$$
(D_{j}(a,p)x)y+x(D_{j+1}(a,p)y)$$
$$
=\{2<q_{j+2}|x>a_{j+2}-
2<a_{j+2}|x>q_{j+2}+
a_{j}(xp_{j})\}y$$
$$
+x\{2<p_{j+3}|y>a_{j+3}-
2<a_{j+3}|y>p_{j+3}+
(q_{j+2}y)a_{j+2}\}$$
$$
=2<q_{j+2}|x>a_{j+2}y-
2<a_{j+2}|x>
q_{j+2}y+
\{a_{j}(xp_{j})\}y$$
$$
+2<p_{j}|y>xa_{j}-
2<a_{j}|y>xp_{j}+
x\{(q_{j+2}y)a_{j+2}\}.$$
However,
we note
$$
\{a_{j}(xp_{j}\}y=
2<a_{j}|y>xp_{j}-
\{y(xp_{j})\}a,$$
$$
=2<a_{j}|y>xp_{j}-
\{2<y|p_{j}>x-
p_{j}(xy)\}a_{j}$$
$$
=2<a_{j}|y>xp_{j}-
2<y|p_{j}>x a_{j}+
\{p_{j}(xy)\}a_{j}$$
and
$$
x\{(q_{j+2}y)a_{j+2}\}=
2<x|a_{j+2}>q_{j+2}y-
a_{j+2}\{(q_{j+2}y)x\}$$
$$
=2<x|a_{j+2}>q_{j+2}y-
a_{j+2}\{2<q_{j+2}|x>y-
(xy)q_{j+2}\}$$
$$
=2<x|a_{j+2}>q_{j+2}y-
2<q_{j+2}|x>a_{j+2}y+
a_{j+2}\{(xy)q_{j+2}\}$$
so that we get
$$
\{D_{j}(a,p)x\}y+
x\{D_{j+1}(a,p)y\}$$
$$
=2<q_{j+2}|x>a_{j+2}y-
2<a_{j+2}|x>q_{j+2}y 
+2<a_{j}|y>xp_{j}-
2<y|p_{j}>xa_{j}$$
$$+\{p_{j}(xy)\}a_{j}
+
2<p_{j}|y>xa_{j}-2<a_{j}|y>xp_{j}
$$
$$
+
2<x|a_{j+2}>q_{j+2}y-
2<q_{j+2}|x>
a_{j+2}y+
a_{j+2}\{(xy)q_{j+2}\}
$$
$$
=\{p_{j}(xy)\}a_{j}
+a_{j+2}\{(xy)q_{j+2}\}=
D_{j+2}(a,p)(xy).
$$
Letting
$j\rightarrow j+1,$
this proves Theorem 3.2.$\square$
\par
For an example of $\Lambda (a)$ in the case of 8 dimensional pseudo octonion
algebra, if $a=(e_{1},e_{2},e_{3}) \in \Sigma $
as in the end of section 2, then we note 
$$
p=(e_{8},e_{8},e_{1}+e_{2}) \in \Lambda (a),
$$
since $e_{1}e_{8}+e_{8}e_{2}=e_{1}+e_{2} $ and $ <e_{8}|e_{1}*e_{2}>, $ etc.

\par
Next, for another proof of Theorem 3.2, 
we will show that
$D_{j}(a,p)$
may be regarded as a local sum of some  element of
$\sigma_{j}(x)\theta_{j}(y)$
which generate the invariant
sub-group of ${\rm Trig}(A)$
as we noted in the previous section.
For this purpose we temporaly
assume that the field $F$ is either real or complex 
field with a formal infinitesimal
variable 
$\varepsilon\in F$
with
$|\varepsilon|<<1.$
We heuristically proceed as follows.
\par
Let $a=(a_{1},a_{2},a_{3}),\ b=(b_{1},b_{2},b_{3})\in \Sigma$
and suppose that
$b\in \Sigma$
is infinitesimally
close to $a\in\Sigma,$
and we write
$$
b_{j}=a_{j}+\varepsilon p_{j}+O(\varepsilon^{2}),\ \ (j=1,2,3)\eqno(3.9)
$$
for some $p_{j}\in A.$
Then,
the condition
$b_{j}b_{j+1}=b_{j+2}$
with
$<b_{j}|b_{j}>=1$
is equivalent to have
$p=(p_{1},p_{2},p_{3})\in \Lambda(a)$
with
$<a_{j}|p_{j}>=0$
and we calculate
$$
\theta_{j}(b)x=
\theta_{j}(a)x+
\varepsilon\{2<p_{j+2}|x>a_{j+1}
+2<a_{j+2}|x>p_{j+1}
-xp_{j}\}+
O(\varepsilon^{2})$$
so that we can write
$$
\sigma_{j}(a)\theta_{j}(b)x=
x+\varepsilon Dx+O(\varepsilon^{2})\eqno(3.10a)$$
for some
$D\in {\rm End} (A)$
defined by
$$
Dx=2<p_{j+2}|x>\sigma_{j}(a)a_{j+1}+
2<a_{j+2}|x>\sigma_{j}(a)p_{j+1}-
\sigma_{j}(a)(xp_{j}).\eqno(3.10b)
$$
We calculate
$$
\sigma_{j}(a)a_{j+1}=
(a_{j+1}a_{j+2})a_{j+1}=
a_{j}a_{j+1}=a_{j+2}$$
$$
\sigma_{j}(a)p_{j+1}=
(p_{j+1}a_{j+2})a_{j+1}=
(p_{j}-a_{j+1}p_{j+2})a_{j+1}
$$
$$
=p_{j}a_{j+1}-<a_{j+1}|a_{j+1}>
p_{j+2}=
p_{j}a_{j+1}-p_{j+2}=
-a_{j}p_{j+1}
$$
so that
$$
Dx=2<p_{j+2}|x>a_{j+2}-2<a_{j+2}|x>a_{j}p_{j+1}$$
$$
-\{2<a_{j+1}|xp_{j}>
a_{j+2}-a_{j}(xp_{j})\}
$$
$$
=2<(p_{j+2}-p_{j}a_{j+1})|x>a_{j+2}-
2<a_{j+2}|x>a_{j}p_{j+1}+
a_{j}(xp_{j})$$
$$
=2<a_{j}p_{j+1}|x>a_{j+2}-
2<a_{j+2}|x>
a_{j}p_{j+1}+
a_{j}(xp_{j})=
D_{j}(a,p)x.$$
Therefore,
Eqs.(2.20)
are rewrite as
$$
\sigma_{j}(a)\theta_{j}(b)=1+\varepsilon D_{j}(a,p)
+O(\varepsilon^{2})\eqno(3.11)$$
for
$D_{j}(a,p)$
defined
by Eq.(3.5)
or (3.7).
\par
Then,
since
$\sigma(a)\theta(b)\in { Trig}(A),$
we have
$$
\sigma_{j}(a)\theta_{j}(b)(xy)=
(\sigma_{j+1}(a)\theta_{j+1}(b)x)y+
x(\sigma_{j+2}(a)\theta_{j+2}(b)y)$$
which implies
$$
\{1+\varepsilon  D_{j}(a,p)+
O(\varepsilon^{2})\}(xy)=
(\{1+\varepsilon D_{j+1}(a,p)+O(\varepsilon^{2})\}x)y+
x\{1+\varepsilon D_{j+2}(a,b)+O(\varepsilon^{2})\}y$$
and hence
$$
D_{j}(a,p)(xy)=
(D_{j+1}(a,p)x)y+
x(D_{j+2}(a,p)y),$$
is given an alternative heuristic approach in the proof of Theorem 3.2.
\par
We note that all relations is given so far as invariant under the action
of $\phi$ defined by
$$
\phi:a_{j}\rightarrow a_{j+1},\ p_{j}\rightarrow p_{j+1},\ q_{j}
\rightarrow q_{j+1},\eqno(3.12)$$
which leads to
$$D_{j}(a,p)\rightarrow D_{j+1}(a,p).$$
\par
Next,
let us consider the relation between
$D_{j}(a,p)$
and
the normal
triality map
$d_{j}(x,y)$
defined by
$$
d_{1}(x,y)=R(y)L(x)-R(x)L(y)\eqno(3.13a)$$
$$
d_{2}(x,y)=L(y)R(x)-L(x)R(y)\eqno(3.13b)$$
while
$d_{3}(x,y)$
for the symmetric composition algebra
is given by
([O,05])
$$
d_{3}(x,y)z=
4\{<x|z>y-
<y|z>x\}.\eqno(3.13c)
$$
If we choose
$p_{3}=0,$
this leads to
$q_{1}=0,\ q_{2}=p_{2}=a_{3}p_{1}$
and $q_{3}=a_{1}p_{2}=-p_{1}a_{2}$
by Eqs.(3.2),(3.3) and (3.4).
For any two arbitary constants
$\alpha,\beta\in F$
with
$\beta\not=0,$
we set then
$$
u={1\over 2\beta}(p_{2}+\alpha a_{2}),\ v=\beta a_{2}.\eqno(3.14)
$$
Then, we can readily verify
the validity of
$$
D_{j}(a,p)=d_{j}(u,v).\eqno(3.15)
$$
Similarly,
if we choose
$p_{1}=0$
or
$p_{2}=0$
and note
Eq.(3.12),
this also
leads to another identity smilar to Eq.(3.15).
Then, the original $D_{j}(a,p)$
can be expressed  as a
linear combinations
of these
$d_{j}(u,v)'s$ for $j=1,2,3$.
\par
For the symmetric composition algebra,
$d_{j}(x,y)$
given by Eqs.(3.13)
satisfy a cubic equation
$$
(d_{j}(x,y))^{3}=
\Delta(x,y)d_{j}(x,y),\ (j=1,2,3)\eqno(3.16a)$$
$$
\Delta(x,y)=
4(<x|y><x|y>-
<x|x><y|y>)\eqno(3.16b)$$
as we may readily verify.
As a matter of fact,
for the specal case of
$j=1,2,$
we have a stronger
relation of
$$
(d_{j}(x,y))^{2}=
\Delta(x,y)Id,\ (j=1,2).\eqno(3.17)$$
Then,
its exponential map can be
computed to be 
$${\rm exp}\{\lambda d_{j}(x,y)\}=$$
$$
1+{sinh(\lambda\sqrt{\Delta(x,y)}\over \sqrt{\Delta(x,y)}}
d_{j}(x,y)+
{1\over \Delta(x,y)}
\{ cosh
(\lambda\sqrt{\Delta(x,y)})
-1\}
(d_{j}(x,y)^{2}.\eqno(3.18)$$
However,
it is harder to be evaluate
${\rm exp}\{\sum^{3}_{j=1}\lambda_{j} D_{j}(a,p)\},\ \ \lambda_{j} \in F$
in a similar way.
\par
\vskip 5mm
{\bf 4\ Automorphism Group of symmetric composition algebras}
\par
\vskip 3mm
In this section, we will discuss the automorphisms of 
the composition symmetric algebra, in particular, for the Cayley algebra,
\par
Let $a\in A$ be a non-trivial idempotent of a 
symmetric composition algebra $A$,
so that
$$
aa=a,\ <a|a>=1.\eqno(4.1)$$
Then,
the triple $(a,a,a)$
is a element of $\sum$ in (2.28) with
$\sigma_{1}(a)=\sigma_{2}(a)=\sigma_{3}(a)$
and
$\theta_{1}(a)=\theta_{2}(a)=\theta_{3}(a).$
In that case,
Theorem 2.5
becomes
\par
\vskip 3mm
{\bf Theorem 4.1}
\par
\vskip 3mm
\it
Let
$a\in A$
be a non-trivial idempotent of $A$,
satisfying Eq.(4.1).
If we set
$$
\sigma(a)=R(a)R(a),\eqno(4.2a)$$
$$
\theta(a)=L(a)L(a),\eqno(4.2b)$$
then they satisfy
\par
(1)
$$\sigma(a)(xy)=(\sigma(a)x)(\sigma(a)y)\eqno(4.3a)$$
\par
(2)
$$\theta(a)(xy)=(\theta(a)x)(\theta(a)y)\eqno(4.3b)$$
\par
(3)
$$
\sigma(a)\theta(a)=
\theta(a)\sigma(a)=Id\eqno(4.3c)$$
\par
(4)
$$
(\sigma(a))^{3}=
(\theta(a))^{3}=Id\eqno(4.3d)$$
\par
(5)
$$
(\sigma(a))^{3}=
(\theta(a))^{3}=Id\eqno(4.3e)$$
\par
(6)
$$
<\sigma(a)x|\sigma(a)y>=
<\theta(a)x|\theta(a)y>=
<x|y>.\eqno(4.3f)$$
Therefore,
both
$\sigma(a)$
and
$\theta(a)$
are automorphisms of $A$ of order $3$.
\par
Moreover
by restricting
$g_{1}=g_{2}=g_{3}(\equiv g)\in { Auto}(A)$
in Eqs.(2.37),
they generate a invariant sub-group of ${ Auto}(A).$
\par
\rm
Especially, from this theorem,
the case of para-Hurwitz algebra with the para-unit $e$ 
is of some interest,
where it can be obtained from Hurwitz algebra
$A^{*}$ i.e., unital composition algebra with the bi-linear product $x*y$ by
$$
xy={\overline {x*y}}=
{\bar y}*{\bar x}\eqno(4.4)$$
with
$e*x=x*e=x,$
and the involution map
$x\rightarrow{\bar x},$
given by
[S.66]
$${\bar x}=2<e|x>e-x.\eqno(4.5)$$
\par
Introducing multiplication operators
$l(x)$
and $r(x)\in {\rm End}A^{*}$
by 
$$
l(x)y=x*y,\ r(x)y=y*x\eqno(4.6)
$$
we find
\par
\vskip 3mm
{\bf Corollary 4.2}
\par
\vskip 3mm
\it
Let $A^{*}$ be a Hurwitz algebra over a field $F$ of charachteristic
$\not= 2$
with the unit element $e$ and the
bi-linear product
$x*y(={\overline{ xy}}).$
For any
$a\in A^{*}$ satisfying
$$
<a|a>=1,\ {\rm and}\ 2<e|a>=-1,\eqno(4.7)$$
we set
$$
\sigma(a)=
l({\bar a})r(a)=
r(a)l({\bar a}).\eqno(4.8)$$
We then have
\par
(i)
$$\overline{\sigma(a)}(x*y)=
(\sigma(a)x)*(\sigma(a)y),\eqno(4.9a)$$
\par
(ii)
$$
\sigma(a)\sigma({\bar a})=Id\eqno(4.9b)$$
\par
(iii)
$$
(\sigma(a))^{3}=Id\eqno(4.9c)$$
\par
(iv)
$$
<\sigma(a)x|\sigma(a)y>=
<x|y>\eqno(4.9d)$$
for any
$x,y\in A^{*}.$
Especially,
$\sigma(a)\in{ Auto}(A^{*}).$
\par
\vskip 3mm
\rm
{\bf Proof}
\par
\vskip 3mm
This is nothing but a special case of Theorem 4.1
by following reason.
We first note that any
$a\in A^{*}$
satisfying Eq.(4.7)
with ${\bar a}=2<a|e>e-a$
satisfies
$a*a={\bar a},$
and hence
$aa=a,$
since any
$x\in A^{*}$
obeys a quadratic relation
$x*{\bar x}=<x|x>e,$
or equivalently
([S,66])
$$
x*x=2<e|x>x-<x|x>e,\eqno(4.10)$$
and
$\sigma(a)=R(a)R(a)$
in $A$.
Then,
Eq.(4.3a)
is rewritten as
$$
\overline{\sigma(a)}(x*y)=
(\sigma(a)x)*(\sigma(a)y).\eqno(4.11)$$
Setting
$y=e,$
this yields
$$\overline{\sigma(a)}x=
(\sigma(a)x)*(\sigma(a)e).$$
However,
$$
\sigma(a)e=l({\bar a})r(a)e={\bar a}*
(e*a)=
{\bar a}*a=<a|a>e=e$$
so that
$\overline{\sigma(a)}x=(\sigma(a)x)*e=
\sigma(a)x$
and then
${\overline \sigma(a)}=\sigma(a).$
\par
Then,
Eq.(4.11)
leads to Eq.(4.9a),
which Eq.(4.8)
gives $\theta(a)=\sigma({\bar a})$
and hence
Eq.(4.9b).This completes the proof.$\square$
\par
\vskip 3mm
{\bf Remark 4.3}
\par
\vskip 3mm
Let $G$ be the group generated by
$\sigma(a)'s$.
Then,
$G$ is a invariant sub-group of ${ Auto}(A^{*}),$
as we see by setting
$g_{1}=g_{2}=g_{3}(\equiv g)\in
{ Auto}(A^{*})$
in Eq.(2.36).
Suppose now that
$A^{*}$ is the Cayley algebra with basis
$e_{0}(\equiv e),e_{1},\cdots ,e_{7}$
satisfying
$$
<e_{\mu}|e_{\nu}>=\delta_{\mu\nu}.\ 
(\mu,\nu=0,1,2,\cdots,7).$$
Then,
expressing
$$
a={1\over 2}(-e+\sum^{7}_{\lambda=1}
\alpha_{\lambda}e_{\lambda})$$
for
$\alpha_{\lambda}\in F,$
the condition
Eq.(4.7)
is equivalent to the
validity of
$$
\sum^{7}_{\lambda=1}
(\alpha_{\lambda})^{2}=3.$$
These imply that 
there is an idempotent element $a$ in the assumption of Theorem 4.1 and Cor.4.2 for the Cayley algebra, since $aa=a$ and $ <a|a>=1$.
Furthermore, 
if the field
$F$ is of charachteristic not $3$,
we can identify
such an element as a point on the
6-sphere
$S^{6}.$
Hence we may be found idempotent  elements
for other algebras by same way.
Moreover for any $b\in A^{*}$ satisfying $<b|b>=1$ and $2<b|e>=-1$,
$c=\sigma (b)a$ preserves
the same condition when
we note
$$<c|c>=<\sigma(b)a|\sigma(b)a>=<a|a>=1,$$
$$
2<c|e>=2<\sigma(b)a|e>=2<a|\sigma({\bar b})e>=2<a|e>=-1.
$$
Therefore $\sigma(b) $ acts upon $S^{6}$,
although it may not be transitive.
However,
for any
$b,c\in S^{6},$
we can find 
$a\in S^{6}$
satisfying
$$
\sigma(a)b=c,$$
provided that we have
$2<b|c>+1\not=0$
with the field $F$ being quadratically closed.
The solution is given by
$$
a={1\over 1+\lambda+\lambda^{2}}
\{\lambda(b+c)-\lambda^{2}e-c*b\}\in S^{6}$$
where
$\lambda\in F$
is a solution of a quadratic equation
$$
1+\lambda+\lambda^{2}=
2<b|c>+1\ (\not=0).$$
\par
If $b$ and $c$ satisfy
$2<b|c>+1=0,$
then we may not find such
$a\in S^{6}.$
However,
in that case,
we choose
$b^{'}\in S^{6}$
satisfying
$2<c|b^{'}>+1\not=0$,
and $2<b|b^{'}>+1\not= 0$ so that there exists
$a_{1},a_{2}\in S^{6}$
satisfying
$$
\sigma(a_{1})b=b^{'}\ {\rm and}\ \sigma(a_{2})b^{'}=c,$$
it leads to
$\sigma(a_{2})\sigma(a_{1})b=c.$
This implies that the group
$G$ acts now
transitively upon
$S^{6}.$
\par
We also note that
$G$ is an invariant sub-goup
of
${ Auto}(A^{*}),$
we have
$t\sigma(a)t^{-1}=\sigma(ta)$
for any
$t\in { Auto}(A^{*}).$
if we choose
$t=\sigma(b).$
This gives
$$
\sigma(b)\sigma(a)=
\sigma(\sigma(b)a)\sigma(b).
$$
\par
The corresponding
 derivation algebra
is obtaind from
Eq.(3.5)
by setting
$a_{1}=a_{2}=a_{3}(\equiv a),\ 
p_{1}=p_{2}=p_{3}(\equiv p),$
and
$q_{1}=q_{2}=q_{3}(\equiv q).$
Then,
$$
D(a,p):=D_{1}(a,p)=D_{2}(a,p)=D_{3}(a,p)$$
is given by
$$
D(a,p)x=
(px)a+
a(xq)=
{\bar a}*(p*x)+
(x*q)*{\bar a}
\eqno(4.12a)$$
or
$$
D(a,p)=
l({\bar a})l(p)+
r({\bar a})r(q).\eqno(4.12b)$$
Note that there are many
constrainted
relations among
$a,p$
and $q$
as we see in the following.
\par
\vskip 3mm
{\bf Lemma 4.4}
\par
\vskip 3mm
\it
Under the assumption as in Corollary 4.2,
we have
\par
(1)
$$
<a|a>=1,\ <p|p>=<q|q>=2<p|q>,\ <p|a>=<q|a>=0,\eqno(4.13a)$$
\par
(2)
$$
2<a|e>=-1,\ <p|e>=<q|e>=0,\ {\bar a}=-e-a,\eqno(4.13b)$$
\par
(3)
$$
q=-p*{\bar a}=-a*p,\ p=-{\bar a}*q=-q*a,\eqno(4.13c)$$
\par
(4)
$$
q*p=<p|p>a,\ p*q=<p|p>{\bar a}.\eqno(4.13d)$$
\par
\vskip 3mm
\rm
{\bf Proof}
\par
\vskip 3mm
From Eq.(3.2),
we calculate
$$
p=ap+pa=
{\bar p}*{\bar a}+
{\bar a}*{\bar p}=
2<{\bar p}|e>{\bar a}+
2<{\bar a}|e>{\bar p}-
2<{\bar a}|{\bar p}>e$$
$$
=2<p|e>{\bar a}-
{\bar p}-0=
2<p|e>{\bar a}-
(2<p|e>e-p)$$
$$
=2<p|e>({\bar a}-e)+p$$
which gives
$<p|e>=0.$
Then,
Eq.(3.3)leads to
$q=ap={\bar p}*{\bar a}=-p*{\bar a}.$
Also,
we calculate
$$
<q|e>=<ap|e>=
<a|pe>=<a|{\bar p}>=
-<a|p>=0$$
and hence
$<q|e>=0.$
Now,
Eq.(3.4)
is rewritten as
$p=qa={\bar a}*{\bar q}=
-{\bar a}*p.$
Moreover,
$$
p*q=
-p*{\bar q}=
({\bar a}*q)*{\bar q}=
{\bar a}<q|q>e=
<q|a>{\bar a}$$
as well as
$$
<q|p>=<q|-{\bar a}*q>=
-<e*q|{\bar a}*q>=
-<e|{\bar a}><q|q>
={1\over 2}(q|q>,$$
$$
<p|a>=
<-q*a|a>=
-<q*a|e*a>=
-<q|e><a|a>=0,$$
$$
<q|a>=
-<a*p|a*e>=
-<a|a><p|e>=0.$$
Finally,
we get
$$
p*q+q*p=
2<p|e>q+
2<q|e>p-
2<p|q>e=
-<p|p>e,
$$
and
these complete the proof.$\square$
\par
\vskip 3mm
{\bf Lemma 4.5}
\par
\vskip 3mm
\it
Under the assumption as in Corollary 4.2,
for any $f,g,x\in A^{*},$ and e is unit element,
we have
\par
(1)
$$
f*(g*x)=
x*(f*g)-
2<f|e>(x*g)+
2<g|e>(f*x)+
2<f|x>g-
2<g|x>f,\eqno(4.14a)$$
\par
(2)
$$
(x*f)*g=
(f*g)*x+
2<f|e>(x*g)-
2<g|e>(f*x)-
2<f|x>g+
2<g|x>f,\eqno(4.14b)
$$
\par
(3)
$$
f*(x*g)=
-x*(f*g)+
2<x|e>(f*g)+
2<f|e>(x*g)-
2<f|x>g,\eqno(4.14c)$$
\par
(4)
$$
(f*x)*g=
-(f*g)*z+
2<x|e>(f*g)+
2<g|e>(f*x)
-2<g|x>f.\eqno(4.14d)
$$
\par
\vskip 3mm
\rm
{\bf Proof}
\par
\vskip 3mm
We calculate
$$
f*(x*g)+
x*(f*g)=
\{2<f|e>x+
2<x|e>f-
2<f|x>e\}*g$$
$$
=2<f|e>x*g+
2<x|e>f*g-
2<f|x>g
$$
and
$$
(f*x)*g+
(f*g)*x=
f*\{2<g|e>x+
2<x|e>g-
2<x|g>e\}
$$
$$
=2<g|e>g*x+
2<x|e>f*g-
2<x|g>f
$$
which
give
Eqs.(4.14c)
and (4.14d).
We proceed
then
$$
f*(g*x)+
f*(x*g)=
f*\{2<x|e>g+
2<g|e>x-
2<x|g>e\}
$$
$$
=2<x|e>f*g+
2<g|e>f*x-
2<x|g>f
$$
and
$$
(f*x)*g+
(x*f)*g=
\{2<x|e>f+
2<f|e>x-
2<f|x>e\}*g$$
$$
=2<x|e>f*g+
2<f|e>x*g-
2<f|x>g
$$
which lead to Eqs.(4.14a)
and (4.14b)
in view of Eqs.(4.14c and d).
$\square$
\par
We can then rewrite Eq.(4.12a) to be
$$
D(a,p)x=
{\bar a}*(p*x)+
(x*q)*{\bar a}$$
$$
=x*({\bar a}*p)-
2<{\bar a}|e>x*p+
2<p|e>({\bar a}*x)+
2<{\bar a}|x>p-
2<p|x>{\bar a}
$$
$$
+(q*{\bar a})*x+
2<q|e>(x*{\bar a})
-2<{\bar a}|e>(q*x)
-
2<q|x>{\bar a}+
2<{\bar a}|x>q$$
$$
=x*({\bar a}*p)+
x*p+
2<{\bar a}|x>p-
2<p|x>{\bar a}$$
$$
+(q*{\bar a})*x+q*x-
2<q|x>{\bar a}+
2<{\bar a}|x>q.$$
\par
Further,
$$
{\bar a}*p=
(-e-a)*p=
-p+q,$$
and
$$
q*{\bar a}=q*(-e-a)=
-q+p$$
so that we find
$$
D(a,p)x=
(x*q+p*x)+
2<{\bar a}|x>(p+q)
-2<(p+q)|x>{\bar a}.\eqno(4.15)
$$
\par
Next,
we would like to relate this derivation to the
standard one
([S,66])
given by
$$
d(x,y)=
l([x,y]^{*})-
r([x,y]^{*})-
3[l(x),r(y)]\eqno(4.16a)
$$
$$=
[l(x),l(y)]+
[r(x),r(y)]+
[l(x),r(y)]\eqno(4.16b)
$$
where
$$[x,y]^{*}\equiv x*y-y*x.\eqno(4.16c)
$$
\par
We can then similarly find
$$
d(f,g)x=
\{-2(f*g)-
(g*f)+
6<g|e>g\}*x$$
$$
+x*\{2(f*g)+
(g*f)-
6<f|e>g\}+
6<f|x>g-6<g|x>f\eqno(4.17)
$$
for any 
$f,g,x\in A^{*}.$
If we choose
$f={\bar a}$
and
$g=p+q$
in Eq.(4.17),
we find
$$
d({\bar a},p+q)=
3D(a,p)\eqno(4.18)
$$
\par
Moreover,
if we set
$$
u=\alpha{\bar a}+\beta(p+q)+\lambda e$$
$$
v=\alpha^{'}{\bar a}+
\beta^{'}(p+q)+\lambda^{'}e$$
for arbitray
$\alpha,\beta,\lambda,\alpha^{'},\beta^{'},\lambda^{'}
\in F$
satisfying
$$\alpha\beta^{'}-\beta\alpha^{'}=1.$$
Then, we can rewrite Eq.(4.18)
as
$$
d(u,v)=3D(a,p).\eqno(4.18)'$$
In this form,
we need not
 worry about constrained relations over $u$ and $v$,
since
$$
<u|e>=
-{1\over 2}\alpha+\lambda,\ <v|e>=
-{1\over 2}\alpha^{'}+\lambda^{'},$$
so that we can regard
$u$
and
$v$
be any generic element of
$A^{*}$.
Therefore,
our derivation
$D(a,p)$
represents the standard formula
for
${ Der}(A^{*}),$
if the underlying field $F$ is of charachteristic
$\not=2,$
and
$\not=3.$
\par
Returning to the discussion of ${ Auto}(A^{*}),$
the most general formula has
been given by
Elduque ([E.00],
who generalized the earlier result of
Jacobson([J.58]):
\par
\vskip 3mm
{\bf Proposition 4.6 ([E.00])}
\par
\vskip 3mm
\it
Let $A^{*}$
be a Cayley-Dickson
algebra
with the unit element
$e$.
\par
(i)\quad
For any $a_{1},\cdots ,a_{r}\ \in A^{*}$ 
and for any positive integer $r$, if it
satisfies
$$
a_{1}*(a_{2}*(\cdots(a_{r-1}*a_{r}))=e=((a_{r}*a_{r-1})
\cdots)*a_{1}.\eqno(4.19).$$
then, we have
$$
{ Auto}(A^{*})=
\{\prod^{r}_{i=1}l(a_{j}),\forall r\in Z,\ \forall a_{i}\in A^{*}\}\eqno(4.20a)$$
$$
=\{
\prod^{r}_{i=1}r(a_{i}),\ \forall r\in Z,\ \forall a_{i}\in A^{*}\}\eqno(4.20b)$$
$$
=\{\prod^{r}_{i=1}
l(a_{i})r(a_{i}),\ \forall r\in Z,\ \forall a_{i}\in A^{*}\}\eqno(4.20c)$$
\par
(ii)\quad
Let the underlying field $F$ be of characteristic
$\not=2$ and
let
$a_{1},a_{2},\cdots,a_{r},$
with
$r=2s=$
even to satisfy
$$
<e|a_{i}>=0,\ a_{1}*(a_{2}*(\cdots (a_{r-1}*a_{r})))=e\eqno(4.21)
$$
instead of Eq.(4.19).
Then,
${ Auto}(A^{*})$
is still expressed as in Eq.(4.20).
Especially,
${ Auto}(A^{*})$ is inner.
\rm
\par
\vskip 3mm
{\bf Remark 4.7}
\par
\vskip 3mm
The relationship between this formula and the
invariant sub-group $G$ of 
${ Auto}(A),$
generated by
$\sigma(a)$'s in Theorem 4.1
is not clear.
In course of its study,
we have found the following rather peculiar 
automorphism
$\sigma\in { Auto}(A^{*})$
satisfying
$\sigma^{2}=2\sigma-1$
from
Eqs.(4.20)
for the case of
$r=3.$
More generally,
we note
\par
\vskip 3mm
{\bf Theorem 4.8}
\par
\vskip 3mm
\it
Let
$A^{*}$ be any algebra over a field $F$ of charachtericstic
$\not=2$,
which needs not be the
Hurwitz algebra.
If $\sigma\in { Auto}(A^{*})$ satisfies
$\sigma^{2}=2\sigma-1,$
then
$d=\sigma-1$
is a derivation of $A^{*}$ satisfying
$dd=0.$
Conversely,
if $d\in { Der}(A^{*})$
satisfies
$dd=0,$
then
$\sigma=d+1$
is an automorphism of 
$A^{*}$,
satisfying
$\sigma^{2}=2\sigma-1.$
\par
\vskip 3mm
\rm
{\bf Proof}
\par
\vskip 3mm
Suppose that
$\sigma\in{\rm Auto}(A^{*}),\sigma \not =1$
satisfies
$\sigma^{2}=2\sigma-1.$
We then calcultate
$$
(2\sigma-1)(x*y)=
\sigma^{2}(x*y)=
\sigma\{(\sigma x)*(\sigma y)\}=
(\sigma^{2}x)*(\sigma^{2}y)$$
$$
=(2\sigma-1)x*(2\sigma-1)y=
4(\sigma x)*(\sigma y)-
2(\sigma x)*y-(2x)(\sigma y)+x*y$$
$$
=4\sigma(x*y)-2(\sigma x)*y-
2x*(\sigma y)+x*y,$$
which is rewritten as
$$
d(x*y)=
(dx)*y+x*(dy)$$
for
$d=\sigma-1,$
which satisfies
$dd=(\sigma-1)^{2}=0,$
also.
\par
Conversely,
if $d\in {\rm Der}(A^{*})$
satisfies
$dd=0,$
then we have
$$
0=dd(x*y)=
d\{dx*y+
x*dy\}=
(ddx)*y+
(dx)*(dy)+
(dx)*
(dy)+
x*(ddy)$$
so that we obtain
$(dx)*(dy)=0.$
We now calculate
$$
\sigma(x*y)=
(d+1)(x*y)=
(dx)*y+
x*(dy)+x*y$$
$$
=(d+1)x*(d+1)y=
\sigma x*\sigma y.$$
\par
Moreover,
$\sigma^{2}=(d+1)^{2}=
2d+1=2\sigma-1.\square$
\par
For any integer
$n=0,\pm 1,\pm 2,\cdots,$
we then find
$$
\sigma^{n}=n(\sigma-1)+1,$$
which leads to
$\sigma^{n}\sigma^{m}=\sigma^{n+m}$
for any integers
$n$ and 
$m$.
If the field $F$ is of charachteristic $p,$
then this yields
$\sigma^{p}=1$
so that
$\sigma$
generates a cyclic group $Z_{p}.$
On the other side,
if $F$ is of characteristic
zero,
then the generated group  is a infinite dimensinal one.
\par
Returning to discussion of Proposition 4.6,
the case of
$r=2$
in Eq.(4.20)
gives only the trivial automorphism
$\sigma=1.$
To see it,
we note that Eq.(4.19)
implies then
$$
a_{1}*a_{2}=e=a_{2}*a_{1}$$
and have
$a_{2}=<a_{2}|a_{2}>{\bar a}_{1}$
with
$<a_{1}|a_{1}>=<a_{2}|a_{2}>=1,$
so that
$$
\sigma=l(a_{1})l(a_{2})=
<a_{2}|a_{2}>l(a_{1})l({\bar a}_{1})=
<a_{2}|a_{2}>
<a_{1}|a_{1}>Id=Id.$$
Similarly,
$r(a_{1})r(a_{2})=1.$
\par
We consider  the case of $r=3,$
which gives
$\sigma\in{ Auto}(A^{*})$
satisfying
$\sigma^{2}=2\sigma-1,$
as we will sketch below.
In this third case,
Eq.(4.19)
becomes
$$
a_{1}*(a_{2}*a_{3})=e=
(a_{3}*a_{2})*a_{1}.\eqno(4.22)
$$
\par
For simplicity, we normarize $a_{j}$'s
to satisfy 
$$
<a_{1}|a_{1}>=<a_{2}|a_{2}>=<a_{3}|a_{3}>=1
\eqno(4.23)$$
by extending the field $F,$
if necessary.
Then
Eq.(4.22)
is equivalemt to the validity of
$$
a_{j}*a_{j+1}=
a_{j+1}*a_{j}=
\overline{a_{j+2}}\eqno(4.24)
$$
with
$a_{j\pm3}=a_{j}.$
Especially,
$(a_{1},a_{2},a_{3})\in \sum$ (see Eq.(2.28)).
Moreover since
$a_{j}*a_{j+1}=a_{j+1}*a_{j},$
the quadratic relation
$$
x*{\bar y}+y*{\bar x}=
2<x|y>e
$$
together with
Eq.(4.23)
gives
$$
a_{j}+\varepsilon_{j+2}a_{j+1}+\varepsilon_{j+1}a_{j+2}=
(2\varepsilon_{j}+<a_{j+1}|a_{j+2}>)e\eqno(4.25a)
$$
where we here set
$$
\varepsilon_{j}=<e|a_{j}>\eqno(4.25b)
$$
for
$j=1,2,3.$
Therefore,
if
$$
{\rm Det}
\left[
\begin{array}{ccc}
1&\varepsilon_{3}&\varepsilon_{2}\\
\varepsilon_{3}&1&\varepsilon_{1}\\
\varepsilon_{2}&\varepsilon_{1}&1
\end{array}
\right]
=1-(\varepsilon_{1}^{2}+\varepsilon_{2}^{2}+\varepsilon_{3}^{2})
+2\varepsilon_{1}\varepsilon_{2}\varepsilon_{3}
\eqno(4.26)
$$
is not zero,
we must have
$a_{j}=\lambda_{j}e\ (j=1,2,3)$
for some
$\lambda_{j}\in F$
satisfying
$\lambda_{1}\lambda_{2}\lambda_{3}=1.$
But then we calcultate
$$
\sigma=
l(a_{1})l(a_{2})l(a_{3})=
r(a_{1})r(a_{2})r(a_{3})=
\prod^{3}_{r=1}l(a_{r})r(a_{r})=1$$
to be a trivial automorphysm.
Hence,
in order to obtain non-trivial
$\sigma,$
we must supose
$$
1-(\varepsilon_{1}^{2}+\varepsilon_{2}^{2}+
\varepsilon_{3}^{2})+
2\varepsilon_{1}\varepsilon_{2}\varepsilon_{3}=0.\eqno(4.27)
$$
\par
Furthermore,
by Eq.(4.25a),
we obtain
$$
(1-\varepsilon_{2}^{2})a_{1}=
(\varepsilon_{1}\varepsilon_{2}-\varepsilon_{3})a_{2}+
\{
2(\varepsilon_{1}-\varepsilon_{2}\varepsilon_{3})
+
<a_{2}|a_{3}>-
\varepsilon_{2}<a_{1}|a_{2}>\}e.$$
Then,
if $\varepsilon_{2}^{2}\not=1,$
this together
with Eq.(4.25a)
for
$j=3,$
we can express both
$a_{1}$
and
$a_{3}$
as some linear conbinations
of $a_{2}$
and
$e$,
which generate a
associative algebra
together with any
$x\in A^{*}$
by Artin's theorem
([S.66]).
We then find again
$$
\sigma=
l(a_{1})l(a_{2})l(a_{3})=
l((a_{1}*a_{2})*a_{3})=l(e)=1$$
to be trivial.
In order to get
a non-trivial
$\sigma,$
we must then
have the relation
$
\varepsilon_{1}^{2}=\varepsilon_{2}^{2}=
\varepsilon_{3}^{2}=1,\ 
\varepsilon_{1}\varepsilon_{2}\varepsilon_{3}=1
$
by Eq(4.26).
\par
{\bf Remark 4.9.}\quad 
Let $\varepsilon_{j} \in F(j=1,2,3)$
be constants to satisfy
$$
\varepsilon_{1}^{2}=\varepsilon_{2}^{2}=
\varepsilon_{3}^{2}=1,\ 
\varepsilon_{1}\varepsilon_{2}\varepsilon_{3}=1.\eqno(4.28)
$$
Also, let
$b_{j}\in A^{*}\ (j=1,2,3)$
be to satisfy
\par
(i)
$$
<b_{i}|e>=0=<b_{i}|b_{j}>\eqno(4.29a)
$$
\par
(ii)
$$
b_{i}b_{j}=0\eqno(4.29b)$$
\par
(iii)
$$
\varepsilon_{1}b_{1}+
\varepsilon_{2}b_{2}+
\varepsilon_{3}b_{3}=0\eqno(4.29c)
$$
for $i,j=1,2,3.$
Then,
$a_{j}\in A^{*}$
defined by
$$
a_{j}=b_{j}+\varepsilon_{j}e\eqno(4.30)$$
satisfy Eq.(4.24).
Moreover, 
$\sigma\in { End}(A^{*})$
given by
$$
\sigma=l(a_{1})l(a_{2})l(a_{3})\eqno(4.31)$$
is an automorphism of $A^{*}$
satisfying $\sigma ^{2}= 2\sigma-1.$
In this case,
we can also
rewrite
$$
\sigma=
1+\varepsilon_{3}l(b_{1})l(b_{2})=
1+\varepsilon_{2}l(b_{3})l(b_{1})=
1+\varepsilon_{1}l(b_{2})l(b_{3}).\eqno(4.32)
$$
We will not go into  detail of this proof.
We can similarly find
also
$$
r(a_{1})r(a_{2})r(a_{3})=
\prod^{3}_{i=1}
l(a_{i})r(a_{i})=\sigma.$$
The next question is whether
non-trivial
$b_{j}'s$
satisfying
Eqs.(4.29)
exist or
not.
The answer is affirmative
if $A^{*}$
is a split
Cayley
algebra,
although we will not go into detail,
however we note that we may find the elements $b_{j}'s$ 
by utilizing of
Eqs (6.11) and (6.12) in ([O.95]).
\par
Finally,
we simply note that
since
$a =(a_{1},a_{2},a_{3})\in\Sigma $ in this case of $r=3$,
we can construct an element of
${ Trig}(A)$
(but not
${ Trig}(A^{*})$)
by
Theorem 2.5.
\par
Next,
let us consider the case of $r=4,$
and note that a special choice of
$a_{4}=e$
reduces the problem 
to that of the case of
$r=3.$
More generally,
the case of
$r=n$
for a integer
$n$ can
be regarded as a special instance of $r=n+1$ 
with the choice of $a_{n+1}=e$.
However the analysis for general case
is complicated 
and we did not succeed in finding some relationship
between Theorem 4.1 and
Proposition 4.6.
In future, we would like to consider for this view point.
\par
\vskip 5mm
{\bf 5\ Examples of triality groups for some nonassociative algebras with involution}
\par
\vskip 3mm
Here,
in this section,
we will give some examples of the triality group for nonassociative algebras with involution
other 
than the symmetric composition algebra.
\par
\vskip 3mm
{\bf Example 5.1\ (Matrix algebra)}
\par
\vskip 3mm
Let
$A=M(n,F)$
be a set consisting all
$n\times n$
matrices over the field $F.$
For any
$x,y\in M(n,F),$
the matrix product which we designate
as
$x*y$
is associative and we write
$$
x*(y*z)=
(x*y)*z:=
x*y*z.\eqno(5.1)$$
Moreover,
for the transpose matrix
$^{t}{x}$
of any
$x\in M(n,F),$
we set
$$
{\bar x}=^{t}{x}.\eqno(5.2)$$
Then
$x\rightarrow {\bar x}$
is a involution map of the resulting algebra
$A^{*}.$
Then,
$A^{*}$ is a unital involutive associative
algebra.
Especially,
it is structurable.
\par
Note that the 
$n\times n$ unit
matrix $e$ is ,
here,
the unit element of $A^{*}.$
We introducea subset of $A^{*}$ by
$$
A_{0}^{*}=
\{x|{\bar x}*x=
x*{\bar x}=e,\ x\in A^{*}\}.\eqno(5.3)
$$
For any three
$a_{j}\in A_{0}^{*},\ (j=1,2,3),$
we introduce
$\sigma_{j}(a)\in {\rm End}A^{*}$
by
$$
\sigma_{j}(a)x:=
a_{j}*x*{\bar a}_{j+1},\ (j=1,2,3)
\eqno(5.4)
$$
where the indices over 
$j$
are defined modulo
$3,$ i.e.,
$a_{j\pm 3}=a_{j}.$
It is easy to see the
validity of
$$
\sigma_{j}(a)\sigma_{j}({\bar a})=1,\eqno(5.5a)$$
$$
\overline{\sigma_{j}(a)}x=
a_{j+1}*x*{\bar a}_{j},\eqno(5.5b)
$$
$$
\overline{\sigma_{j}(a)}(x*y)=
(\sigma_{j+1}(a)x)*
(\sigma_{j+2}(a)y).\eqno(5.5c)
$$
Introducing the conjugate algebra
$A$
of $A^{*}$
with the bi-linear product
$xy$
by
$$
xy=\overline{x*y}=
{\bar y}*{\bar x}.\eqno(5.6)
$$
then, these are rewritten as
$$
\sigma_{j}(a)x=
a_{j+1}(a_{j}x)=
(x{\bar a}_{j+1}){\bar a}_{j}\eqno(5.7a)
$$
$$
\sigma_{j}(a)(xy)=
\{(\sigma_{j+1}(a)x)\}\{
(\sigma_{j+2}(a)y)\}.\eqno(5.7b)
$$
We note that
$A$
is 
no longer associative
but para-associative
with
the para-associative law of
$$
{\bar z}(xy)=
(yz){\bar x}.\eqno(5.8)
$$
From Eqs.(5.5a)
and (5.7b),
we find
$$
\sigma(a)=
(\sigma_{1}(a),\sigma_{2}(a),\sigma_{3}(a))\in
{\rm Trig}(A).\eqno(5.9)
$$
Furthermore,
Eq.(5.7a)
gives
$$
\sigma_{j}(a)=
L(a_{j+1})L(a_{j})=
R({\bar a}_{j})R({\bar a}_{j+1}),\eqno(5.10)
$$
which has the same structures as Eq.(2.29)
for the symmetric
composition algebra.
\par
For the corresponding to local triality  case,
let
$$
p=(p_{1},p_{2},p_{3})\in (A^{*})^{3},\ {\bar p}_{j}
=-p_{j}\ (j=1,2,3)\eqno(5.11)
$$
and define
$d_{j}(p)\in {\rm End}\ A^{*}$
by
$$
d_{j}(p)x=
p_{j}*x-
x*p_{j+1}.\eqno(5.12)$$
We then have
$$
\overline{d_{j}(p)}x=
x*{\bar p}_{j}
-{\bar p}_{j+1}*x=
p_{j+1}*x-
x*p_{j}\eqno(5.13a)$$
$$
\overline{d_{j}(p)}(x*y)=
(d_{j+1}(p)x)*y+
x*(d_{j+2}(p)y),\eqno(5.13b)$$
and hence
$$
d_{j}(p)(xy)=
(d_{j+1}(p)x)y+
x(d_{J+2}(p)y).\eqno(5.14)
$$
These imply that $d(p)=(d_{1}(p),d_{2}(p),d_{3}(p) )\in s\circ Lrt (A).$
\par
This example implies that 
for the matrix algebra, we have  a generalization of  the corresponding
$$
a*x*a^{-1}  \ ( \ automorphism ) \Longleftrightarrow p*x-x*p\  (\ derivation)
$$
where $a$ is  orthogonal and $p$ is alternative,
by means of well-known Cayley transformation, that is,
$a=(Id -p)*(Id +p)^{-1}$ (by assuming well-defined).
\par 
However,
we will not go into detail
\par
It seems that this local triality relation represented by
Eq. (5. 14) may be regarded as a  version for matrix algebras which correspondense with the well known
" the principle of triality" for the Cayley algebra( [S.66]),
and also similarly that the relation (3.6) of Theorem 3.2 may be regarded as a local triality relation (the principle of triality) for the symmetric composition algebra.
\vskip 3mm
{\bf Example 5.2\ (Para-Zorn Matrix Algebra)}
\par
\vskip 3mm
Let $B$ be an algebra over a field $F$ with a bi-linear
form
$(\cdot|\cdot)$
and let
$$
A=\left(\begin{array}{cc}
F&B\\
B&F\end{array}
\right)\eqno(5.15)
$$
be the Zorn's vector matrix with its
generic elemnt
$X\in A$
given
by
$$
X=
\left(
\begin{array}{cc}
\alpha&x\\
y&\beta
\end{array}
\right)\eqno(5.16)
$$
for
$\alpha,\beta\in F$
and
$x,y\in B.$
We introduce a bi-linear product in 
$A$ by
$$
X_{1}X_{2}=
\left(\begin{array}{cc}
\alpha_{1}&x_{1}\\
y_{1}&\beta_{1}
\end{array}\right)
\left(\begin{array}{cc}
\alpha_{2}&x_{2}\\
y_{2}&\beta_{2}\end{array}
\right)=
\left(\begin{array}{ll}
\beta_{1}\beta_{2}+(y_{1}|x_{2}),&
\alpha_{1}x_{2}+\beta_{2}x_{1}+ky_{1}y_{2}\\
\alpha_{2}y_{1}+\beta_{1}y_{2}+kx_{1}x_{2},&
\alpha_{1}\alpha_{2}+(x_{1}|y_{2})
\end{array}
\right)\eqno(5.17)
$$
for $k\in F.$
\par
If $B$ is involutive with the involution map
$x\rightarrow {\bar x}$
 satisfying
$$
\overline{xy}={\bar y}{\bar x},\ 
({\bar x}|{\bar y})=
(y|x),\eqno(5.18)
$$
then $A$ becomes involutive also by
$$
X\rightarrow {\bar X}=
\left(
\begin{array}{cc}
\beta&{\bar x}\\
{\bar y}&\alpha
\end{array}\right),\eqno(5.19)
$$
and its
conjugate algebra $A^{*}$ with the
bi-linear product
$$
X_{1}*X_{2}=\overline{X_{1}X_{2}}=
{\bar X_{2}}{\bar X_{1}}\eqno(5.20)
$$
will present the standard
Zorn's vector matrix algebra.
\par
Also we note that if the involution of $B$ is identity map,
then  $B$ is commutative, moreover assume that 
if $B$ is the exceptional Jordan algebra with $dim\ B=27$,
then 
for $k=2$, $A^{*}$ is a structurable algebra ([A.78], [K-O.15] ).
\par
If the involution of $B$ is anti-commutative (i.e., ${\bar x}=-x$ )
and $B$ is provided the property  equipped with some conditions which satisfy $x(yz)= (x|y)z-(x|z)y $
and  $(x|yz)=(y|zx)=(z|xy)$, then
it is known that for $k=1$
$A^{*}$ is an alternative  algebra ([K-O.15]).
\par
On the other hand, by utilizing of these Zorn's vector matrices,
we have considerd a construction of simple
Lie algebras or superalgebras
from concept of  triple systems with a ternary product $(X_{1}X_{2}X_{3})$ in $A^{*}$,
however we will not go into the details (for example, to see 
[K-O.15] and references therein ).
\par 
Next for any non-zero $\lambda \in F,$
we intoduce
$\rho_{j}(\lambda)\in{ End}\ A\ (j=1,2,3)$
by
$$
\rho_{1}(\lambda)X=
\left(\begin{array}{cc}
\lambda\alpha,&\lambda x\\
\lambda^{-1}y,&
\lambda^{-1}\beta
\end{array}\right)
\eqno(5.21a)
$$
$$
\rho_{2}(\lambda)X=
\left(
\begin{array}{cc}
\lambda \alpha,&\lambda^{-1}x\\
\lambda y,&
\lambda^{-1}\beta
\end{array}
\right)
\eqno(5.21b)
$$
$$
\rho_{3}(\lambda)X=
\left(
\begin{array}{cc}
\lambda^{-2}\alpha,&
x\\
y,&\lambda^{2}\beta
\end{array}
\right).\eqno(5.21c)
$$
We can then show
\par
\vskip 3mm
{\bf Theorem 5.3}
\par
\vskip 3mm
\it
Under the assumption as in Ex.5.2, we have
\par
(i)\quad
$$
\rho_{j}(\lambda)(XY)=
(\rho_{j+1}(\lambda)X)
(\rho_{j+2}(\lambda)Y)\eqno(5.22a)
$$
for any
$X,Y\in A$ and
for any
$j=1,2,3$
with
$\rho_{j\pm 3}(\lambda)=\rho_{j}(\lambda).$
\par
(ii)
$$
\rho_{j}(\mu)\rho_{j}(\nu)=
\rho_{j}(\mu\nu),\ \rho_{j}(1)=1\eqno(5.22b)
$$
\par
(iii)
$$
\rho_{j}(\mu)\rho_{k}(\nu)=
\rho_{k}(\nu)\rho_{j}(\mu)\eqno(5.22c)$$
\par
(iv)
$$
\rho_{1}(\lambda)\rho_{2}(\lambda)
\rho_{3}(\lambda)=1\eqno(5.22d)$$
\par
(v)\quad
If $B$ is involutive,
then
$$
\overline{\rho_{j}(\lambda)}=
\rho_{3-j}(\lambda^{-1})\eqno(5.22e)
$$
for any
$j,k=1,2,3$
and for any non-zero
$\mu,\nu\in F.$
\par
\vskip 3mm
\rm
{\bf Proof}
\par
\vskip 3mm
Since the calculations are straightforward,
we will verify only a few of these statements.
First
$$
(\rho_{1}(\lambda)X_{1})
(\rho_{2}(\lambda)X_{2})=
\left(
\begin{array}{cc}
\lambda\alpha_{1},&\lambda x_{1}\\
\lambda^{-1}y_{1},&
\lambda^{-1}\beta_{1}
\end{array}
\right)
\left(
\begin{array}{cc}
\lambda\alpha_{1},&
\lambda^{-1}x_{2}\\
\lambda y_{2},&
\lambda^{-1}\beta_{2}
\end{array}
\right)
$$
$$
=
\left(
\begin{array}{cc}
\lambda^{-2}\beta_{1}\beta_{2}+
\lambda^{-2}(y_{1}|x_{2}),&
\alpha_{1}x_{2}+\beta_{2}x_{1}+
ky_{1}y_{2}\\
\alpha_{2}y_{1}+
\beta_{1}y_{2}+
kx_{1}x_{2},&
\lambda^{2}\alpha_{1}\alpha_{2}+
\lambda^{2}(x_{1}|y_{2})
\end{array}
\right)
=\rho_{3}(\lambda)(X_{1}X_{2})$$
which proves the case of
$j=3$
for
Eq.(5.22a).
Similarly,
we calcultate
$$
(\rho_{2}(\lambda)X_{1})
(\rho_{3}(\lambda)X_{2})=
\left(
\begin{array}{cc}
\lambda\alpha_{1},&\lambda^{-1}x_{1}\\
\lambda y_{1},&
\lambda^{-1}\beta_{1}
\end{array}
\right)
\left(
\begin{array}{cc}
\lambda^{-2}\alpha_{2},&
x_{2}\\
y_{2},&
\lambda^{2}\beta_{2}
\end{array}
\right)
$$
$$
=
\left(
\begin{array}{cc}
\lambda\beta_{1}\beta_{2}+
\lambda(y_{1}|x_{2}),&
\lambda\alpha_{1}x_{2}+\lambda\beta_{2}x_{1}+
\lambda ky_{1}y_{2}\\
\lambda^{-1}\alpha_{2}y_{1}+
\lambda^{-1}\beta_{1}y_{2}+
\lambda^{-1}kx_{1}x_{2},&
\lambda^{-1}\alpha_{1}\alpha_{2}+
\lambda^{-1}(x_{1}|y_{2})
\end{array}
\right)
=
\rho_{1}(\lambda)(X_{1}X_{2}),$$
for the case of $j=1.$\par
We also note for example
$$
\rho_{2}(\lambda)\rho_{3}(\lambda) X=
\rho_{2}(\lambda)
\left(
\begin{array}{cc}
\lambda^{-2}\alpha,&x\\
y,&\lambda^{2}\beta
\end{array}
\right)=
\left(
\begin{array}{cc}
\lambda(\lambda^{-2}\alpha),&
\lambda^{-1}x\\
\lambda y,&
\lambda^{-1}(\lambda^{2}\beta)
\end{array}
\right)
=\rho_{1}(\lambda^{-1})X$$
so that
$$
\rho_{1}(\lambda)\rho_{2}(\lambda)\rho_{3}(\lambda)X=
\rho_{1}(\lambda)\rho_{1}(\lambda^{-1})X=X.$$
Finally, if $B$ is involutive,
then we compute
$$
\overline{\rho_{1}(\lambda)}{\bar X}=
\overline{\rho_{1}(\lambda)X}=
\left(
\begin{array}{cc}
\lambda^{-1}\beta,&
\lambda{\bar x}\\
\lambda^{-1}{\bar y},&
\lambda \alpha
\end{array}
\right).
$$
On the other side,
$$
\rho_{2}(\lambda^{-1}){\bar X}=
\rho_{2}(\lambda^{-1})
\left(
\begin{array}{cc}
\beta&{\bar x}\\
{\bar y}&\alpha
\end{array}
\right)=
\left(
\begin{array}{cc}
\lambda^{-1}\beta,&
\lambda{\bar x}\\
\lambda^{-1}{\bar y},&
\lambda\alpha
\end{array}
\right)
$$
and
hence
$\overline{\rho_{1}(\lambda)}=
\rho_{2}(\lambda^{-1}).$
We can similarly verify other relations.$\square$
\par
\vskip 3mm
{\bf Remark 5.4}
\par
\vskip 3mm
Let us set
$$
e=\left(
\begin{array}{cc}
1&0\\
0&1
\end{array}
\right),\ 
g(\lambda)=
\left(
\begin{array}{cc}
\lambda&
0\\
0&{1\over \lambda}
\end{array}
\right),\ 
h(\lambda)=
\left(
\begin{array}{cc}
{1\over \lambda}&0\\
0&\lambda
\end{array}
\right).\eqno(5.23)
$$
Then, we have
$$
eX=Xe=
\left(
\begin{array}{cc}
\beta,&x\\
y,&\alpha
\end{array}
\right)\eqno(5.24)
$$
for
$X=
\left(
\begin{array}{cc}
\alpha&x\\
y&\beta
\end{array}
\right),$
and 
$eg(\lambda)=h(\lambda),
g(\lambda)h(\lambda)=e, h(\lambda)e=g(\lambda),$
as in the property
 of $a_{j+2}=a_{j}a_{j+1}$
in section 2,
provided $<e|e>:=det\ e =1, <g(\lambda)|g(\lambda)>:=det\ g(\lambda)=1, <h(\lambda)|h(\lambda)>:=det\ h(\lambda) =1$.
\par
Moreover,
we can express
$\rho_{j}(\lambda)$
as
$$
\rho_{1}(\lambda)=L(e)L(g(\lambda))=
R(e)R(h(\lambda))\eqno(5.25a)$$
$$
\rho_{2}(\lambda)=
L(h(\lambda))L(e)=
R(g(\lambda))R(e)\eqno(5.25b)
$$
$$
\rho_{3}(\lambda)=
L(g(\lambda))L(h(\lambda))=
R(h(\lambda))
R(g(\lambda))\eqno(5.25c)$$
which have structures similar to Eqs.(2.29)
for the symmetric composition algebra,again.
\par
Hence these imply that
$\rho (\lambda)=(\rho_{1}(\lambda),\rho_{2}(\lambda),\rho_{3}(\lambda))\in Trig\  (A).$
\par
Next,
let
$\pi\in { End}(A)$
be given by
$$
\pi:
\left(
\begin{array}{cc}
\alpha&x\\
y&\beta
\end{array}
\right)
\rightarrow
\left(
\begin{array}{cc}
\alpha&y\\
x&\beta
\end{array}
\right)\eqno(5.26)
$$
which satisfies
\par
(i)
$$
\pi(XY)=
(\pi X)(\pi Y)\eqno(5.27a)
$$
\par
(ii)
$$
\pi^{2}=1\eqno(5.27b)$$
\par
(iii)
$$
\pi\rho_{j}(\lambda)\pi^{-1}=\rho _{3-j}(\lambda)=
{\bar \rho_{j}}(\lambda^{-1}).\eqno(5.27c)
$$
Hence the group generated by $\rho (\lambda)=(\rho_{1}(\lambda),
\rho_{2}(\lambda),\rho_{3}(\lambda) )$ and $ \pi_{0}=(\pi,\pi,\pi)$
will make a sub-group of $Trig\ (A)$. 
\par
Moreover suppose that $\xi, \ \eta\in Epi (A)$
satisfy 
\par
(iv)
$$
\xi(xy)=(\eta x)(\eta y)\eqno(5.28a)
$$
\par
(v)
$$
\eta(xy)=(\xi x)(\xi y)\eqno(5.28b)
$$
for any
$x,y\in B.$
We call the pair
$(\xi,\eta)$
be a double automorphism of $B$.
We then find

\par
\vskip 3mm
{\bf Proposition 5.5}
\par
\vskip 3mm
\it
Under the assumption as in Theorem 5.3,
if a double automorphism
$(\xi,\eta)$
of $B$ satisfy
$$
(\xi x|\eta y)=(x|y),\eqno(5.29)
$$
then $P\in { End}A$
given by
$$
P:X\rightarrow
\left(
\begin{array}{cc}
\alpha,&\xi x\\
\eta y,&\beta
\end{array}
\right)\eqno(5.30a)
$$
is an automorphism of the para-Zorn algebra $A$,
i.e.,
we have
$$
P(XY)=(PX)(PY).\eqno(5.30b)
$$
\rm
{\bf Proof.}
\par
However,
we will not go into its proof,
since it is straightforward. $\square$
\par
\vskip 3mm
For simple example of a double automorphism, we have $(\xi,\eta)=(\omega,\omega^{2})$, where $\omega^{3}=1$.
\vskip 3mm
Finally,
let us consider an example of  the corresponding local
triality Lie reated triple in the para-Zorn marix algebra $A$.
We will  introduce
$s_{j}\in { End}\ A$
for
$j=1,2,3$
by
$$
s_{1}X=
\left(
\begin{array}{cc}
\alpha&x\\
-y,&-\beta
\end{array}
\right)\eqno(5.31a)
$$
$$
s_{2}X=
\left(
\begin{array}{cc}
\alpha&-x\\
y&-\beta
\end{array}
\right)\eqno(5.31b)
$$
$$
s_{3}X=
\left(
\begin{array}{cc}
-2\alpha&0\\
0&2\beta
\end{array}
\right).\eqno(5.31c)
$$
Then by the straightfoward calculations, we may obtain the following.
\par
\vskip 3mm
{\bf Proposition 5.6}
\par
\vskip 3mm
\it
Under the assumpsion as in above, we have
\par
(i)
$$
s=(s_{1},s_{2},s_{3})\in s\circ Lrt(A),i.e.,$$
$$
s_{j}(XY)=(s_{j+1}X)Y+
X(s_{j+2}Y)\eqno(5.32a)$$
\par
(ii)
$$
s_{1}+s_{2}+s_{3}=0\eqno(5,32b)$$
\par
(iii)
$$
[s_{j},s_{k}]=0,\eqno(5.32c)$$
\par
(iv)\quad If $B$ is involutive,
then
$$
{\bar s}_{j}=-s_{3-j}\eqno(5.32d)
$$
for any
$j,k=1,2,3,$
where the indices over
$j$
are again defined modulo $3.$
\par
\rm
{\bf Proof.}
\par
Eqs(5.32) are clear by straightfoward calculations.
From ${\bar s_{j}}(X)=\overline{s_{j}({\bar X})}=-s_{3-j}(X),$
it is east to show that
${\bar s_{j}}=-s_{3-j} \ (j=0,1,2)$. $\square$
\par

\vskip 3mm
In ending this section,
we see that three examples considered in this note of the symmetric composition algebra, 
the matrix algebra and
the para-Zorn algebra admit the "Triality group" of form given as in Eq.(2.29).
Hence these suggest that
there may exist
other class of algebras admitting a construction of 
triality groups in this way.
\par
On the other hand, for triple systems equipped with a ternary product $(xyz)$ ( for example, to see [O.95],  [K-O.00], [K-M-O.10]),
it seems that local and global triality relations may be made as
 same concept,
but it will be considered in future.
\par 
\vskip 3mm
{\bf Appendix  (Simple examples)}
\par
In this appendix,
we give simple couple examples.
\par 
{\bf Example A }\ \ Let {\bf C} be the complex number with usual product $x*y$.
If we define for $a =(\alpha,\beta,\gamma)$ $\in \  {\bf C}^{3}$,
and $|\alpha|=| \beta|=| \gamma |=1$,
$$\sigma_{1}(a)x=\alpha *x *\beta^{-1}$$
$$\sigma_{2}(a)x=\beta *x *\gamma^{-1}$$
$$\sigma_{3}(a)x=\gamma * x *\alpha^{-1}$$
where $\sigma_{j\pm 3}=\sigma_{j} (j=0,1,2)$,
then we have
$$ \sigma_{j}(a)(xy)=
(\sigma_{j+1}(a)x) (\sigma_{j+2}(a)y)
$$  
with respect to  new product $xy$  defined by
$xy =\overline{x * y}$,
\par 
\noindent
where ${\bar x}$ denotes the conjugation of $x$.
\par
If we introduce for $a=(\alpha.\beta,\gamma) \in\ (Im\ {\bf C})^{3}$
$$d_{1}(a)x=\alpha *x-x*\beta, 
d_{2}(a)x=\beta *x-x*\gamma,
d_{3}(a)x=\gamma *x-x*\alpha,$$
then we get 
$$d_{j}(a)(xy)=(d_{j+1}(a)x)y +x(d_{j+2}(a)y).
$$
{\bf Example B}\ \ 
Let  {\bf H} be the quaternion algebra satisfying $i^{2}=j^{2}~=k^{2}=-1 , i*j=-j*i=k$ with 
basis $\{1, i,j,k\}$.
If we define for $a=(i,j,k) \in ( Im\  {\bf H})^{3},$
$$\sigma_{1}(a)x=i*x*j^{-1}$$
$$\sigma_{2}(a)x=j*x*k^{-1}$$
$$\sigma_{3}(a)x=k*x*i^{-1}$$
where the product $x*y$ is usual product of {\bf H},
and $\sigma _{j\pm 3}=\sigma _{j}$,
\par
\noindent
then we have 
$$
\sigma_{j}(a)(xy)=(\sigma_{j+1}(a)x)(\sigma_{j+2}(a)y)
$$
w.r.t. new product $xy=\overline{x*y}$,
where ${\bar x}$ denotes the involutive conjugation of $x \in {\bf H}$. 
\par
If we introduce
$d_{1}(a)x=i*x-x*j, d_{2}(a)x=j*x-x*k,d_{3}(a)x=k*x-x*i,$
then we obtain
$$
d_{l}(a)(xy)=(d_{l+1}(a)x)y +x(d_{l+2}(a)y),\ \ l=l\pm 3
$$
w.r.t new product $xy$ defined by $xy=\overline{x*y}.$
\par
{\bf Concluding Remark }
\par
In this above ${\bf H}$, ${\bf H}$  has a structure of a symmetric composition algebra with respect to new product defined by $xy=\overline{x*y}$, Thus special examples of $\Sigma$ in Eq. (2.28) are $(i,j,-k),\ (i,-j,k),\ (-i,j,k)$.
Indeed, from fact that ${\bf H}$ is the symmetric composition algebra with respect to the product $xy$ and 
$ij=\overline{i*j}=-k$ etc, we obtain
 a simple  special  case of Theorem 2.5.
Furthermore, putting $a={1\over 2}(-1+\sqrt{3} \ i)$,
then we have $a*{\bar a}=1$, $aa=a$, $<a|a>=1$,
 because $aa={1 \over 4}\overline{(-1+\sqrt {3} \ i)*(-1+\sqrt {3}\ i)}=a$
and so $a$ is an idempotent element,
hence  it satisfies the assumpsion of Theorem 4.1.
That is, if we set  $\sigma (a)=R(a)R(a)$, then $\sigma (a)$ is
an automorphism, i.e.,
$$
\sigma (a)(xy)=(\sigma(a)x)(\sigma(a)y).
$$
{\bf Acknowledgments}
\par
The one of authors (N.Kamiya) would like to express his gratitude to
Prof. A. Elduque for some correspondences.
\par
 \vskip 5mm
{\bf References}
\par
\rm
\vskip 3mm
\noindent
[A.78]:
Allison,B.N.;
"A class of non-associative algebras with involution
containing
a class of Jordan algebra"\ 
Math,Ann 
{\bf 237},
(1978)
133-156
\par
\noindent
[A-F.93]:
Allison,B.N.,Faulkner;J.R.;
"Non-associative coefficient algebras for
Steinberg
unitary Lie algebras"\ 
J.Algebras {\bf 161}
(1993)133-158
\par
\noindent
[E.97]:
Elduque,A.;
"Symmetric Composition Algebra"\ 
J.Algebra
{\bf 196}(1997)
282-300
\par
\noindent
[E.00]:
Elduque,A.;
"On triality and automorphism and derivation of
composition algebra"\ 
Linear algebra and 
its applications
{\bf 314}(2000)
49-74
\par
\noindent
[G.62]
Gell-Mann,M.;
"Symmetry of Baryons and Mesons, Phys. Rev., vol. 125
(1962) 1067-1084.

\par
\noindent
[J.58]:
Jacobson,N.:
"Composition algebras and thier automorphisms"\ 
Rend. Circ.
Mat. Palermo 7 (1958),55-80
\par
\noindent
[K-M-O.10]
Kamiya,N., Mondoc,D.. Okubo.S.,;
"Structure theory of  (-1,-1)Freudenthal-Kantor triple systems",
Bull.Austr.Math.Soc., {\bf 81} (2010) 132-155. 
\par
\noindent
[K-M-P-T.98]:
Knus,M.A.,Merkurjev,A.S.,Post,M.,Tignal,J.P.,;
"The Book of Involution"\ 
American Math.Soc.Coll.Pub.
{\bf 44}
Providence
(1998)
\par
\noindent
[K-O.00] Kamiya,N., Okubo,S., ;
" On $\delta$- Lie super triple systems associated with $(\varepsilon,\delta)$ Freudenthal-Kantor triple systems", 
Proc. Edinburgh Math.Soc., {\bf 43} (2000) 243-260
\par
\noindent
[K-O.14]
Kamiya,N., Okubo,S., 
"Triality of structurable and pre-structurabe algebras" J.Alg.,
{\bf 416}(2014) 58-88 
\par
\noindent
[K-O.15]:
Kamiya,N. Okubo,S.;
"Algebras satisfying triality and 
$S_{4}$
Symmetry"
 (2015) Arxiv.1503.00614,
Algebras, groups and geometries, vol {\bf 33}, nu.1, (2016)1-92.
\par
\noindent
[O.95]
Okubo,S. ;
"Introduction to octonion andother non-associativealgebras to physics" Cambridge Univ, press, Cambridge(1995)
\par
\noindent
[O.05]:
Okubo,S.;
"Symmetric triality relations and 
structurable algebra"\ 
Linear Algebras and
its Applications,
{\bf 396}
(2005)
189-222
\par
\noindent
[O-O.81]:
Okubo,S.,Osborn,J.M.;
"Algebras with non-degenerate associative symmetric 
bi-linear form
permitting compositions"\ 
Comm.Algebra
{\bf 9}
(1981)
(I) 1233-1261,
(II) 2015-2073
\par
\noindent
[S.66]:Schafer,R.D.;
"An Introduction to Non-associate algebra"\ 
Academic press.
N.Y. and London
(1966)
\end{document}